\documentclass[11pt]{article}

\usepackage{amsmath,amsfonts,amssymb}
\usepackage{bbm}
\usepackage{cases}
\usepackage{color}
\usepackage[hidelinks]{hyperref}

\def\|{\vert\vert}

\def\bR {{\mathbb R}}

\def\<{{\langle}}
\def\>{{\rangle}}

\newtheorem{theorem}{Theorem}[section]
\newtheorem{lemma}[theorem]{Lemma}

\newtheorem{proposition}[theorem]{Proposition}

\newtheorem{remark}[theorem]{Remark}
\newtheorem{example}[theorem]{Example}
\newtheorem{definition}[theorem]{Definition}

\begin{document}

\title{\vspace*{-1.5cm}
Global well-posedness to stochastic reaction-diffusion equations on the real line $\bR$ with superlinear drifts driven by multiplicative space-time white noise}

\author{Shijie Shang$^{1}$ and
Tusheng Zhang$^{2}$}
%\footnotetext[1]{\, EPF-Lausanne}
%\footnotetext[2]{\, University of Utah}
%\footnotetext[3]{\, School of Mathematics, University of Science and Technology of China, Hefei, China.}
\footnotetext[1]{\, School of Mathematics, University of Science and Technology of China, Hefei, China. Email: sjshang@ustc.edu.cn}
%\footnotetext[4]{\, School of Mathematics, University of Manchester, Oxford Road, Manchester
%M13 9PL, England, U.K. Email: tusheng.zhang@manchester.ac.uk}
\footnotetext[2]{\, Department of Mathematics, University of Manchester, Oxford Road, Manchester
M13 9PL, England, U.K. Email: tusheng.zhang@manchester.ac.uk}
\maketitle

\begin{abstract}
Consider the stochastic reaction-diffusion equation with logarithmic nonlinearity driven by space-time white noise:
\begin{numcases}{}
 \mathrm{d}u(t,x) = \frac{1}{2}\Delta u(t,x) \,\mathrm{d}t+ b(u(t,x)) \,\mathrm{d}t \nonumber\\
~~~~~~~~~~~~~ + \sigma(u(t,x))  \,W(\mathrm{d}t,\mathrm{d}x),  \   t>0, x\in I , \nonumber\\
 u(0,x)=u_0(x), \quad  x\in I .\nonumber
\end{numcases}
When $I$ is a compact interval, say $I=[0,1]$, the well-posedness of the above equation was established in \cite{DKZ} (Ann. Prob. 47:1,2019). The case where $I=\bR$ was left open. The essential obstacle  is caused  by the explosion of the supremum norm of the solution, $\sup_{x\in\bR}|u(t,x)|=\infty$, making the usual truncation procedure invalid. In this paper, we prove that there exists a unique global solution to the stochastic reaction-diffusion equation on the whole real line $\mathbb{bR}$ with logarithmic nonlinearity.  Because of the nature of the nonlinearity, to get the uniqueness, we are forced to work with the first order moment of the solutions on the space $C_{tem}(\bR)$ with a specially designed norm
$$\sup_{t\leq T, x\in\mathbb{R}}\left(|u(t,x)|e^{-\lambda |x|e^{\beta t}}\right),$$
where, unlike the usual norm in $C_{tem}(\bR)$, the exponent also depends on time $t$ in a particular way.
 Our approach depends heavily on  the new, precise lower order moment estimates of the stochastic convolution and a new type of Gronwall's inequalities we obtained, which are of  interest on their own right.
\end{abstract}


\medskip


\tableofcontents

\section{Introduction}

%In this paper, we study the equation with a superlinear growth term driven by space-time white noise. This paper is divided into two parts.
%
%In the first part of this paper, we study the following stochastic equation,
%
%
%the heat equation with a logarithmic nonlinear term driven by additive space-time white noise,
%\begin{numcases}{}
% du(t,x) = \frac{1}{2}\Delta u(t,x) dt+ u(t,x)\log|u(t,x)| dt \nonumber\\
%~~~~~~~~~~~~~ + \sigma  W(dt,dx),  \   t>0, x\in \mathbb{R} , \nonumber\\
%\label{1.a} u(0,x)=u_0(x), \quad  x\in \mathbb{R} .
%\end{numcases}
%The coefficient $\sigma$ is a constant. $W$ is space-time white noise on $\mathbb{R}_{+}\times \mathbb{R}$ defined on some filtrated probability space $(\Omega, {\cal F}, \{{\cal F}_t\}_{t\geq 0}, \mathbb{P})$.
%
%{\color{red}
%First, we consider the additive noise case, i.e. $\sigma$ is a constant. The multiplicative noise case is TBA.}
%

In this paper, we study the stochastic reaction-diffusion equation on the whole line $ \bR$ driven by multiplicative space-time white noise given as follows:
\begin{numcases}{}
 \mathrm{d}u(t,x) = \frac{1}{2}\Delta u(t,x) \,\mathrm{d}t+ b(u(t,x)) \,\mathrm{d}t \nonumber\\
~~~~~~~~~~~~~ + \sigma(u(t,x))  \,W(\mathrm{d}t,\mathrm{d}x),  \   t>0, x\in \mathbb{R} , \nonumber\\
\label{1.a} u(0,x)=u_0(x), \quad  x\in \mathbb{R} .
\end{numcases}
The coefficients $b, \sigma$ are two deterministic measurable functions from $\mathbb{R}$ to from $\mathbb{R}$, $W$ is a space-time white noise on $\mathbb{R}_{+}\times \mathbb{R}$ defined on some filtrated probability space $(\Omega, {\cal F}, \{{\cal F}_t\}_{t\geq 0}, \mathbb{P})$.
\vskip 0.3cm
 There exist numerous work in the literature on stochastic reaction-diffusion equations driven by space-time white noise covering a wide range of topics.
% ranging from well-posedness, existence and uniqueness of invariant measures, large deviation to  intermittency, Malliavin calculus etc.
 We refer the reader to \cite{DPZ}, \cite{cerrai-1}, \cite{K} and references therein. The majority of the work are focused on  stochastic reaction-diffusion equations defined on finite intervals (i.e., the space variable belongs to a fixed finite interval) instead of the whole real line $\mathbb{R}$, partly due to the essential difficulties brought by the non-compactness of the whole space. We like to mention some relevant existing work on the well-posedness of the stochastic reaction-diffusion equations on the real line. In the  early paper \cite{S}, the author obtained the existence and uniqueness of solutions of stochastic reaction-diffusion equations on the real line under the Lipschitz conditions of the coefficients. Pathwise uniqueness were established in \cite{MP} and \cite{MPS} for stochastic reaction-diffusion equations on the real line with H\"older continuous coefficients.
\vskip 0.4cm
It is well known that the equation (\ref{1.a}) admits a unique global solution when the coefficients fulfill the usual Lipschitz condition, in particular, being of linear growth.  We are concerned here with the well-posedness of the stochastic reaction-diffusion equation  (\ref{1.a}) with superlinear  drift. Several papers in the literature discuss stochastic partial differential equations with locally Lipschitz coefficients that have polynomial growth and/or satisfy certain monotonicity conditions (see \cite{cerrai-1,DMP,LR}, for instance). The typical example of such a coefficient  is $b(u) = -u^3$, which has the effect of ``pulling the solution back toward the origin.'' In the joint paper \cite{DKZ} with Dalang and Khoshnevisan by the second  named author ,  stochastic reaction-diffusion equations(SRDEs) on finite intervals were considered and  it was proved  that
if the coefficients are locally Lipschitz and of  $(|z|\log|z|)$-growth, then  the SRDEs is globally well-posed. Unfortunately, the methods in \cite{DKZ} are not valid  for SRDEs on the whole line $\mathbb{R}$ because typically the supremum norm of the solution explodes, i.e., $|u(t)|_{\infty}=\sup_{x\in \bR}|u(t,x)|=\infty$  . The global well-posedness on the whole line $\mathbb{R}$ under the logarithmic nonlinearity  was left open.
\vskip 0.4cm
The goal of this article is to fill in this gap. More precisely, we  prove that
if the drift $b$ is locally Log-Lipschitz and if $|b(z)|=O(|z|\log|z|)$, then the stochastic reaction-diffusion equation \eqref{1.a} is globally well-posed. The precise statements are given in the next section. Because of the nature of the nonlinearity, we are forced to work with the first order moment of the solutions on the space $C_{tem}(\bR)$ with a specially designed norm
$$\sup_{t\leq T, x\in\mathbb{R}}\left(|u(t,x)|e^{-\lambda |x|e^{\beta t}}\right),$$
where, unlike the usual norm in $C_{tem}(\bR)$, the exponent also depends on time $t$ in a particular way. We need to establish some new, precise (lower order) moment estimates of stochastic convolution on the real line  and hence obtain some a priori estimates  of the solution. We  like to stress that it is harder to get  precise lower order moment estimate than high order  for stochastic convolutions. To obtain the pathwise uniqueness, one of the  difficulties is that we are not able to apply the usual localization procedure as in the literature  because the usual uniform norm of the solution on the real line explodes. To overcome the difficulty, we  provide a new type of Gronwall's inequalities, which is of independent interest.

\vskip 0.4cm
 Now we describe the content of the paper in more details.   In Section 2, we present the framework for the stochastic reaction-diffusion equations driven by  space-time white noise and state our main results. In Section 3, we will prove two  Gronwall-type inequalities and  obtain some estimates associated with the heat kernel of the Laplacian operator. In Section 4, we establish new lower order moment estimates of the stochastic convolution with respect to the space-time white noise and obtain a priori estimates of the solutions of the stochastic reaction-diffusion equations. In Section 5, we approximate the coefficients $b$ and $\sigma$ by smooth functions and establish the tightness of the laws of the solutions of the corresponding approximating equations. As a consequence, we obtain the existence of weak solution ( in the probabilistic sense). Section 6 is devoted to the proof of the pathwise uniqueness of the stochastic reaction-diffusion equation under the local Log-Lipschitz conditions of the coefficients.

\section{Statement of main results}\label{S:2}
\setcounter{equation}{0}

Let us recall the following definition.

\begin{definition}
A random field solution to equation (\ref{1.a}) is a jointly measurable and adapted space-time process $u:=\{u(t,x): (t,x)\in\mathbb{R}_{+}\times\mathbb{R}\}$ such that for every $(t,x)\in \mathbb{R}_{+}\times \mathbb{R}$,
\begin{align}\label{definition solution}
  u(t,x)= & P_t u_0(x)+\int_0^t\int_{\mathbb{R}} p_{t-s}(x,y) b(u(s,y)) \,\mathrm{d}s\mathrm{d}y \nonumber\\
  & + \int_0^t\int_{\mathbb{R}} p_{t-s}\sigma(u(s,y)) \,W(\mathrm{d}s,\mathrm{d}y), \quad \mathbb{P}-a.s.,
\end{align}
where $p_{t}(x,y):=\frac{1}{\sqrt{2\pi t}} e^{-\frac{(x-y)^2}{2t}}$, and $\{P_t\}_{t\geq 0}$ is the corresponding heat semigroup on $\mathbb{R}$.
\end{definition}
\begin{remark}
The above mild form is equivalent to the weak (in the sense of partial differential equations) formulation of the stochastic reaction-diffusion equations. We refer readers to \cite{WA} for details.
\end{remark}
We also recall the so-called $C_{tem}$ space defined by
\begin{align*}
  C_{tem}:=\left\{f\in C(\mathbb{R}): \sup_{x\in\mathbb{R}}|f(x)|e^{-\lambda |x|}<\infty \text{ for any } \lambda>0 \right\},
\end{align*}
and endow it with the metric defined by
\begin{align*}
  d(f,g):=\sum_{n=1}^{\infty}\frac{1}{2^n}\min\left\{1, \sup_{x\in\mathbb{R}}|f(x)-g(x)|e^{-\frac{1}{n} |x|}\right\} ,
\end{align*}
for any $f,g\in C_{tem}$.
Then $f_n\rightarrow f$ in $C_{tem}$ iff $\sup_{x\in\mathbb{R}}|f_n(x)-f(x)|e^{-\lambda |x|}\rightarrow 0$ as $n\rightarrow\infty$ for any $\lambda>0$, and $(C_{tem}, d)$ is a Polish space.

%Let us introduce a new space $L^{1}_{tem}$ which is defined by
%\[
%  L^1_{tem}:=\left\{f\in \mathcal{B}(\mathbb{R}): \int_{\mathbb{R}}|f(x)|e^{-\lambda |x|}dx<\infty \text{ for any } \lambda>0 \right\},
%\]
%where $\mathcal{B}(\mathbb{R})$ is the set of all Borel measurable functions on $\mathbb{R}$. The metric on $L^1_{tem}$ is defined by
%\begin{align*}
%  d_{L^1}(f,g):=\sum_{n=1}^{\infty}\frac{1}{2^n}\min\left\{1, \int_{\mathbb{R}}|f(x)-g(x)|e^{-\frac{1}{n} |x|}dx\right\} .
%\end{align*}

\vskip 0.3cm

Next we introduce the following conditions of nonlinear term $b$. Set $\log_{+}(u):=\log_{+}(1\vee u)$ for any $u\geq 0$.
\begin{itemize}
	\item [(H1)] $b$ is continuous, and there exist two nonnegative constants $c_1$ and $c_2$ such that for any $u\in\mathbb{R}$,
 \begin{align}\label{210124.2000}
  |b(u)|\leq c_1 |u|\log_{+}|u| + c_2.
 \end{align}
	\item [(H2)] There exist nonnegative constants $c_3, c_4, c_5$, such that for any $u,v\in\mathbb{R}$,
 \begin{align}
  |b(u)-b(v)|\leq c_3|u-v|\log_{+}\frac{1}{|u-v|}+  c_4 \log_{+}(|u|\vee|v|)|u-v| +  c_5|u-v|.
 \end{align}
\end{itemize}

Note that  condition (H2) implies condition (H1). A typical example of function $b$ that satisfies (H2) is given below.
\vskip 0.4cm
\begin{example}\label{A.2}
The function $x\mapsto x\log|x|$ satisfies the local log-Lipschitz condition (H2), more precisely, for any $x,y\in\mathbb{R}$,
\begin{align}\label{05161157}
   \left|x\log|x|-y\log|y|\right|\leq |x-y|\log\frac{1}{|x-y|} + [\log_{+}(|x|\vee|y|)  +1 +\log 2]|x-y|.
\end{align}
\end{example}
\noindent {\bf Proof}.
Without loss of generality, we may assume $|y|\leq |x|$. We divide the proof into two cases.

Case 1: $|x| \leq 1$. We have
\begin{align}
  \left|x\log|x|-y\log|y|\right|\leq & |x-y|\left|\log|x|\right|+ |y|\left|\log|x|-\log|y|\right| \nonumber\\
  \leq & |x-y|\log\frac{1}{|x|} + |x-y| \nonumber\\
  \leq & |x-y|\left(\log\frac{1}{|x-y|} +\log 2\right) +|x-y| ,
\end{align}
since $|x-y|\leq 2|x|$ leads to
\begin{align}
  \log\frac{1}{|x|}\leq \log\frac{1}{|x-y|} + \log 2.
\end{align}

%Case 2. $1\leq |x|\leq M$. We have
%\begin{align}
%  \left|x\log|x|-y\log|y|\right|\leq & |x-y|\left|\log|x|\right|+ |y|\left|\log|x|-\log|y|\right| \nonumber\\
%  \leq & |x-y|\times \log M  +|x-y| .
%\end{align}
Case 2. $|x|\geq 1$. We have
\begin{align}
  \left|x\log|x|-y\log|y|\right|\leq & |x-y|\left|\log|x|\right|+ |y|\left|\log|x|-\log|y|\right| \nonumber\\
  \leq & |x-y|\times \log_{+}(|x|\vee|y|)  +|x-y| .
\end{align}
Combining these two cases together yields (\ref{05161157}).
$\blacksquare$

\vskip 0.5cm

Now we can state the main results of this paper.

\begin{theorem}\label{thm1}
Assume  $u_0 \in C_{tem}$ and that (H1) is satisfied. If $\sigma$ is bounded and continuous, then there exists a  weak ( in the probabilistic sense) solution to the stochastic reaction-diffusion equation (\ref{1.a}) with sample paths a.s. in $C(\mathbb{R}_{+}, C_{tem})$.
\end{theorem}

%\begin{theorem}\label{thm2}
%Assume that condition (H1) is satisfied. If $\sigma$ is continuous and satisfies the globally linear growth condition, $u_0 \in C_{tem}\cap L^1_{tem}$, then there exists a solution to (\ref{1.a}), and the sample paths of the solution are a.s. in $C(\mathbb{R}_{+}, C_{tem}\cap L^1_{tem})$.
%\end{theorem}

\vskip 0.6cm

\begin{theorem}\label{thm3}
	Assume  $u_0 \in C_{tem}$ and that (H2) is satisfied. If $\sigma$ is bounded and Lipchitz, then the pathwise uniqueness holds for solutions of (\ref{1.a}) in $C(\mathbb{R}_{+}, C_{tem})$. Hence there exists a unique strong solution to (\ref{1.a}) in $C(\mathbb{R}_{+}, C_{tem})$.

%the pathwise uniqueness of solutions holds in $C(\mathbb{R}_{+}, C_{tem})$ for equation (\ref{1.a}).
\end{theorem}

\section{Preliminaries}
\setcounter{equation}{0}

In this section, we will provide two  Gronwall-type inequalities which play an important role in this paper. Moreover, we also present some estimates associated with the heat kernel of the Laplacian operator which will be used in our analysis later.
\vskip 0.3cm
Lemma \ref{A.1} is a slight modification of Theorem 3.1 in \cite{W}, and is proved in Lemma 7.2 of \cite{SZ2}. We give a short proof  here for completeness. Set $\log_{+}(r):= \log(r\vee 1)$.

\begin{lemma}\label{A.1}
Let $X, a, c_1, c_2$ be nonnegative functions on $\mathbb{R}_{+}$, $M$ an increasing function with $M(0)\geq 1$. Moreover, suppose that  $c_1, c_2$ be integrable on finite time intervals. Assume that for any $t\geq 0$,
\begin{align}\label{5.3}
  X(t)+a(t)\leq M(t)+\int_0^t c_1(s)X(s)\,\mathrm{d}s +\int_0^t c_2(s) X(s)\log_{+} X(s) \,\mathrm{d}s ,
\end{align}
and the above integral is finite. Then for any $t\geq 0$,
\begin{align}\label{5.4}
  X(t)+a(t)\leq M(t)^{\exp(C_2(t))}\exp\left(\exp(C_2(t))\int_0^t c_1(s)\exp(-C_2(s))\,\mathrm{d}s\right) ,
\end{align}
where $C_2(t):=\int_0^t c_2(s)\,\mathrm{d}s $.
\end{lemma}
%\begin{remark}\label{05161918}
%  If the last term of (\ref{5.3}) is $\int_0^t c_2(s) X(s)\log_{+} X(s) ds$, where $\log_{+}(r):= \log(r\vee 1)$, then the estimate (\ref{5.4}) still holds for the same reason.
%\end{remark}
\noindent {\bf Proof}. Fix any $T>0$. Let
\[
Y(t):= M(T)+\int_0^t c_1(s)X(s)\,\mathrm{d}s +\int_0^t c_2(s) X(s)\log_{+} X(s) \,\mathrm{d}s, \quad t\in[0,T] .
\]
We see that $Y$ is almost surely differentiable on $[0,T]$, $Y(t)\geq 1$ and
\[
X(t)+a(t)\leq Y(t), \quad \forall\, t\in [0,T] .
\]
This leads to
\begin{align}\label{5.5}
  Y^{\prime}(t)= & c_1(t)X(t)+ c_2(t) X(t)\log_{+} X(t) \nonumber\\
  \leq &  c_1(t)Y(t)+ c_2(t) Y(t)\log_{+} Y(t) \nonumber\\
  = &  c_1(t)Y(t)+ c_2(t) Y(t)\log Y(t) .
\end{align}
Thus,
\begin{align}\label{5.6}
  \left(\log Y\right)^{\prime}(t)\leq c_1(t)+ c_2(t)\log Y(t) .
\end{align}
Solving this ordinary differential inequality, we get for any $t\in [0,T]$,
\begin{align}\label{5.7}
  \log Y(t) \leq \exp(C_2(t))\left[\log M(T) + \int_0^t c_1(s)\exp(-C_2(s))\,\mathrm{d}s\right] .
\end{align}
Therefore, we obtain
\begin{align}\label{5.8}
   & X(T)+a(T) \leq  Y(T) \nonumber\\
\leq & M(T)^{\exp(C_2(T))}\exp\left(\exp(C_2(T))\int_0^T c_1(s)\exp(-C_2(s))\,\mathrm{d}s\right) .
\end{align}
By the arbitrariness of $T$, (\ref{5.4}) is deduced.
$\blacksquare$

\begin{lemma}\label{A.3}
Let $Y(t)$ be a nonnegative function on $\mathbb{R}_+$. Let $c_1$ and $c_2$ be non-negative, increasing functions on $\mathbb{R}_+$. Let $\varepsilon\in [0,1)$ be a constant and $c_3: \mathbb{R}_{+}\times(\varepsilon,1)\longmapsto\mathbb{R}_{+}$ be a function that is increasing with respect the first variable. Suppose that for any $\theta\in (\varepsilon,1)$, the following integral inequality holds
\begin{align}\label{210217.2130}
  Y(t)\leq c_1(t)\int_0^t Y(s)\,\mathrm{d}s + c_2(t)\int_0^t Y(s)\log_{+}\frac{1}{Y(s)}\,\mathrm{d}s + c_3(t,\theta)\int_0^t Y(s)^{\theta} \,\mathrm{d}s ,\quad \forall\, t\geq 0 .
\end{align}
If for any $t>0$,
\begin{align}\label{210219.1709}
	\limsup_{\theta\rightarrow 1-} \, (1-\theta) c_3(t,\theta)  <\infty ,
\end{align}
then $Y(t)= 0$ for any $t\geq 0$. In particular, if $c_3(t,\theta) \leq \frac{c(t)}{1-\theta}$ and $c$ is an increasing function with respect to $t$, then (\ref{210219.1709}) holds.
\end{lemma}
\noindent {\bf Proof}.
It suffices to prove that for any $T>0$, $Y(\cdot)=0$ on $[0,T]$. In order to prove this,  let
\begin{align}
	\delta_T := \limsup_{\theta\rightarrow 1-} \, (1-\theta) c_3(T,\theta), \nonumber\\
	T^* : = \min \left\{ T, \frac{1}{3\delta_T}, \frac{e}{3 c_2(T)} \right\} .
\end{align}
%\begin{align}
%	\delta_T :=\max \left\{\limsup_{\theta\rightarrow 1-} \, (1-\theta) c_3(T,\theta) ,\frac{2c_2(T)}{e} \right\}
%\end{align}
%let
%%$\delta^*:=\max\left\{t: c_2(t)t\leq \frac{e}{2}\right\}$,
%$T^*=\frac{1}{\delta^T}$.

Step 1. We first prove $Y(t)=0$ for any $t\in[0,T^*]$. Since
\begin{align}
	\sup_{x\geq 0}\left(x\log_{+}\frac{1}{x}\right)=\frac{1}{e} ,
\end{align}
we have
\begin{align}
	  Y(t)\leq & c_1(t)\int_0^t Y(s)\,\mathrm{d}s + \frac{c_2(t)}{1-\theta}\int_0^t Y(s)^{\theta}Y(s)^{1-\theta}\log_{+}\frac{1}{Y(s)^{1-\theta}}\,\mathrm{d}s \nonumber\\
	  & + c_3(t,\theta)\int_0^t Y(s)^{\theta} \,\mathrm{d}s  \nonumber\\
	  \leq & c_1(t)\int_0^t Y(s)\,\mathrm{d}s + \left[ \frac{c_2(t)}{e(1-\theta)} + c_3(t,\theta)\right] \int_0^t Y(s)^{\theta}\,\mathrm{d}s ,\quad \forall\, t\geq 0 .
\end{align}
For $t\in[0,T]$, let
\begin{align}
	\Phi(t) := c_1(T) \int_0^t Y(s) \,\mathrm{d}s + \left[ \frac{c_2(T)}{e(1-\theta)} + c_3(T,\theta)\right] \int_0^t Y(s)^{\theta}\,\mathrm{d}s .
\end{align}
Then $Y(t)\leq \Phi(t)$ for any $t\in[0,T]$. Thus,
\begin{align}
	\frac{\mathrm{d}}{\mathrm{d}t}\Phi(t) = & c_1(T) Y(t) + \left[ \frac{c_2(T)}{e(1-\theta)} + c_3(T,\theta)\right] Y(t)^{\theta} \nonumber\\
	\leq &  c_1(T)\Phi(t) + \left[\frac{c_2(T)}{e(1-\theta)} + c_3(T,\theta)\right] \Phi(t)^{\theta} .
\end{align}
Without loss of generality, we can assume that $\Phi(t)>0$ for any $t\in(0,T]$, otherwise we can take the zero time to be $\min\{t: \Phi(t)>0\}$.
Multiplying $(1-\theta)\Phi(t)^{-\theta}$ on both sides of the above inequality yields
\begin{align}
	\frac{\mathrm{d}}{\mathrm{d}t}\left(\Phi(t)^{1-\theta}\right) \leq (1-\theta) c_1(T)\Phi(t)^{1-\theta} +   \left[ \frac{c_2(T)}{e} + c_3(T,\theta)(1-\theta)\right].
\end{align}
Solving the above inequality, we obtain
\begin{align}
	\Phi(t)^{1-\theta}\leq \left[ \frac{c_2(T)}{e} + c_3(T,\theta)(1-\theta)\right] \int_0^t e^{(1-\theta)c_1(T)(t-s)} \,\mathrm{d}s .
\end{align}
Hence
\begin{align}
	Y(t)\leq & \Phi(t) \leq \left\{ \left[\frac{c_2(T)T^*}{e} + (1-\theta) c_3(T,\theta)T^* \right] e^{(1-\theta)c_1(T) T^*}  \right\}^{\frac{1}{1-\theta}} \nonumber\\
	\leq & e^{c_1(T)T^*}  \left\{  \frac{c_2(T)T^*}{e} + (1-\theta) c_3(T,\theta)T^* \right\}^{\frac{1}{1-\theta}} ,
\end{align}
for any $t\in[0,T^*]$. Letting $\theta\rightarrow 1$ in the above inequality and in view of the definition of $T^*$,  we can see that
\begin{align}
	Y(t)=0,\quad \forall\, t\in[0,T^*].
\end{align}

Step 2.
From (\ref{210217.2130}) it follows that
\begin{align}
	Y(t+T^*)\leq & c_1(T)\int_0^{t+T^*} Y(s)\,\mathrm{d}s + c_2(T)\int_0^{t+T^*} Y(s)\log_{+}\frac{1}{Y(s)}\,\mathrm{d}s \nonumber\\
	& + c_3(T,\theta)\int_0^{t+T^*} Y(s)^{\theta} \,\mathrm{d}s  \nonumber\\
	\leq & c_1(T)\int_{T^*}^{t+T^*} Y(s)\,\mathrm{d}s + c_2(T)\int_{T^*}^{t+T^*} Y(s)\log_{+}\frac{1}{Y(s)}\,\mathrm{d}s \nonumber\\
	& + c_3(T,\theta)\int_{T^*}^{t+T^*} Y(s)^{\theta} \,\mathrm{d}s   \nonumber\\
	\leq & c_1(T)\int_{0}^{t} Y(s+T^*)\,\mathrm{d}s + c_2(T)\int_{0}^{t} Y(s+T^*)\log_{+}\frac{1}{Y(s)}\,\mathrm{d}s \nonumber\\
	& + c_3(T,\theta)\int_{0}^{t} Y(s+T^*)^{\theta} \,\mathrm{d}s ,\quad  t\in [0, T-T^*] .
\end{align}
By Step 1, we see that
\begin{align}
	Y(t+T^*) = 0, \quad t\in[0, T^*\wedge(T-T^*)] ,
\end{align}
that is $Y(t)=0$ for any $t\in[0, T\wedge 2T^* ]$. Repeating this argument, we see that $Y(t)=0$ for any $t\in [0,T]$. The arbitrariness of $T$ leads  $Y(t)=0$ for any $t>0$.
$\blacksquare$

\vskip 0.3cm

Recall that $p_{t}(x,y):=\frac{1}{\sqrt{2\pi t}} e^{-\frac{(x-y)^2}{2t}}$. Now we state some estimates of the heat kernel $p_t(x,y)$ used in this paper.
%Before we state Lemma 3.3 about estimates of heat kernels,
%We notice the following three estimates, which can be seen by straightforward calculations.
For any $x\in\mathbb{R}$ and $t>0$,
\begin{align}
\label{05132049.1} \int_{\mathbb{R}} p_{t}(x,y)e^{\eta |y|} \,\mathrm{d}y \leq & 2e^{\frac{\eta^2 t}{2}}e^{\eta |x|}, \quad \forall\ \eta\in\mathbb{R} ,  \\
  \label{2102121.2100} \int_{\mathbb{R}} p_{t}(x,y)^2 e^{\eta |y|} \,\mathrm{d}y \leq & \frac{1}{\sqrt{\pi t}}e^{\frac{\eta^2 t}{4}}e^{\eta |x|}, \quad \forall\ \eta\in\mathbb{R},  \\
\label{05161718} \int_{\mathbb{R}}p_{t}(x,y) e^{\eta |y|} \eta |y|  \,\mathrm{d}y \leq & e^{\frac{\eta^2 t}{2}}e^{\eta |x|}\eta|x| + 2e^{\frac{\eta^2 t}{2}}\left(\eta^2 t+ \eta\sqrt{\frac{t}{2\pi}}\right) e^{\eta |x|} ,
 \quad \forall\ \eta >0 .
\end{align}
The above three estimates can be obtained through straightforward calculations. Moreover, the following
lemma is needed.
%
%In the sequel, we will use the following bounds, which can be seen by straightforward calculations.
%Note that for any $\eta, t>0$ and $x\in\mathbb{R}$,

\begin{lemma}\label{A.4}
%Let $p_s(z)=\frac{1}{\sqrt{2\pi s}}e^{-\frac{z^2}{2s}}$, $p_s(x,z)=p_s(x-z)$, $s>0$ and $x, z\in\mathbb{R}$.
%For any $x,y,z\in\mathbb{R}$ and $t>0$, we have the following estimates.
%Let $\eta$ be a positive constant.  We have the following estimates.
The following estimates of the heat kernel $p_{t}(x,y)$ hold.
\begin{itemize}
  \item [(i)] For any $x,y\in\mathbb{R}$, $\theta\in[0,1]$, $0< s\leq t$,
  \begin{align}
  |p_t(x,y)-p_s(x,y)|\leq \frac{(2\sqrt{2})^\theta |t-s|^{\theta}}{s^{\theta}}\big(p_s(x,y)+p_t(x,y)+p_{2t}(x,y)\big) .
\end{align}
  \item [(ii)] For any $x,y\in\mathbb{R}$ and $t>0$,
  \begin{align}
    \int_{\mathbb{R}}|p_t(x,z)-p_t(y,z)|\,\mathrm{d}z\leq \sqrt{\frac{2}{\pi}}\times \frac{|x-y|}{\sqrt{t}}.
  \end{align}
  \item [(iii)] For any $x,y\in\mathbb{R}$ and $\eta, t>0$,
  \begin{align}
    \int_{\mathbb{R}}|p_t(x,z)-p_t(y,z)|e^{\eta |z|}\,\mathrm{d}z\leq 2\sqrt{2}\times \frac{|x-y|}{\sqrt{t}}\times  e^{\eta^2 t}\times e^{\eta (|x|+|x-y|)}.
  \end{align}
  \item [(iv)] For any $x,y\in\mathbb{R}$ and $\eta, t>0$,
  \begin{align}
    & \int_{\mathbb{R}}|p_t(x,z)-p_t(y,z)|e^{\eta |z|}\eta |z| \,\mathrm{d}z \nonumber\\
    \leq & \frac{\sqrt{2}|x-y|}{\sqrt{t}} \times \bigg[ e^{\eta^2 t}\times e^{\eta (|x|+|x-y|)}\eta (|x|+|x-y|) \nonumber\\
    & + 2e^{ \eta^2 t }\bigg(2 \eta^2 t+ \eta\sqrt{\frac{t}{\pi}}\bigg) e^{\eta (|x|+|x-y|)}\bigg].
  \end{align}
  \item [(v)] For any $x,y\in\mathbb{R}$ and $0<s\leq t$,
  \begin{align}\label{210222.1639}
  	\int_0^{s} \int_{\mathbb{R}} |p_{t-r}(x,z) -p_{s-r}(y,z)|^2 \,\mathrm{d}r\mathrm{d}z \leq \frac{\sqrt{2}-1}{\sqrt{\pi}} |t-s|^{\frac{1}{2}} + \frac{2}{\sqrt{\pi}}|x-y| .
  \end{align}
\end{itemize}
\end{lemma}
\noindent {\bf Proof}.
Proof of (i).  On the one hand, by the mean value theorem, there exists some $\xi\in[s,t]$ such that
\begin{align}\label{05161649.1}
  & |p_t(x,y)-p_s(x,y)| \nonumber\\
  \leq & \frac{1}{\sqrt{2\pi t}}\left|e^{-\frac{(x-y)^2}{2t}}-e^{-\frac{(x-y)^2}{2s}}\right| + \left|\frac{1}{\sqrt{2\pi t}}-\frac{1}{\sqrt{2\pi s}}\right|e^{-\frac{(x-y)^2}{2s}} \nonumber\\
  \leq & \frac{1}{\sqrt{2\pi t}} e^{-\frac{(x-y)^2}{2\xi}}\times \frac{(x-y)^2}{2\xi^2}\times|t-s|+ \frac{1}{\sqrt{2\pi}}\times \frac{|t-s|}{2s^{3/2}}e^{-\frac{(x-y)^2}{2s}} \nonumber\\
\leq & \frac{2}{\sqrt{2\pi t}} e^{-\frac{(x-y)^2}{2\xi}}\times e^{\frac{(x-y)^2}{4\xi}}\times \frac{|t-s|}{s} + \frac{|t-s|}{2s}p_s(x,y) \nonumber\\
  \leq & 2\sqrt{2}\times \frac{|t-s|}{s}p_{2t}(x,y)  + \frac{|t-s|}{2s}p_s(x,y) \nonumber\\
  \leq & \frac{2\sqrt{2}|t-s|}{s}\big(p_s(x,y)+p_t(x,y)+p_{2t}(x,y)\big) ,
\end{align}
where we have used $z\leq e^z$ for any $z\geq 0$. On the other hand,
\begin{align}\label{05161649.2}
  |p_t(x,y)-p_s(x,y)|\leq & p_s(x,y)+ p_t(x,y) + p_{2t}(x,y) .
\end{align}
Combining (\ref{05161649.1}) and (\ref{05161649.2}) together yields (i).

Proof of (ii). Obviously,
\begin{align}\label{05161703}
  |p_t(x,z)-p_t(y,z)| = & \frac{1}{\sqrt{2\pi t}}\bigg|\int_0^1 d e^{-\frac{|x-z+\rho (y-x)|^2}{2t}}\bigg| \nonumber\\
  \leq & \frac{|x-y|}{\sqrt{2\pi t}} \times \int_0^1 e^{-\frac{|x-z+\rho (y-x)|^2}{2t}} \times \frac{|x-z+\rho (y-x)|}{t} \,\mathrm{d}\rho  .
%  \leq & \sqrt{\frac{2}{s}}|x-y|\times \int_0^1 \frac{1}{\sqrt{2\pi\times 2s}} e^{-\frac{|x+\rho (y-x)-z|^2}{4s}} d\rho,
\end{align}
Due to the Fubini theorem,
\begin{align}
	& \int_{\mathbb{R}}|p_t(x,z)-p_t(y,z)|\,\mathrm{d}z \nonumber\\
	\leq & \int_0^1 \frac{|x-y|}{t} \int_{\mathbb{R}}\frac{1}{\sqrt{2\pi t}} |x-z+\rho (y-x)| e^{-\frac{|x-z+\rho (y-x)|^2}{2t}} \,\mathrm{d}z \mathrm{d}\rho \nonumber\\
	\leq & \int_0^1 \frac{|x-y|}{t} \sqrt{\frac{2t}{\pi}} \,\mathrm{d}\rho = \sqrt{\frac{2}{\pi t}}|x-y| .
\end{align}
%
%
%where we have used $x\leq e^{\frac{x^2}{4}}$ for any $x\geq 0$.
%Now (ii) follows easily from (\ref{05161703}) and the Fubini theorem.

Proof of (iii).
%By (\ref{05132049.1}) and the Fubini theorem, we have
%\begin{align}
%      & \int_{\mathbb{R}}|p_s(x,z)-p_s(y,z)|e^{\eta |z|}dz \nonumber\\
%      \leq &  \sqrt{\frac{2}{s}}|x-y| \int_0^1 \int_{\mathbb{R}} \frac{1}{\sqrt{2\pi\times 2s}} e^{-\frac{|x+\rho (y-x)-z|^2}{4s}} e^{\eta|z|} dz d\rho \nonumber\\
%      \leq &  \sqrt{\frac{2}{s}}|x-y| \times \int_0^1 2 e^{\eta^2 s}\times e^{\eta |x+\rho (y-x)|} d\rho  \nonumber\\
%      \leq & 2\sqrt{2}\times \frac{|x-y|}{\sqrt{s}}\times e^{\eta^2 s}\times e^{\eta (|x|+|x-y|)} .
%\end{align}
Note that $z\leq e^{\frac{z^2}{4}}$ for any $z\geq 0$. By (\ref{05161703}), we have
\begin{align}\label{05161703.1}
  |p_t(x,z)-p_t(y,z)|  \leq  \sqrt{\frac{2}{t}}|x-y|\times \int_0^1 \frac{1}{\sqrt{2\pi\times 2t}} e^{-\frac{|x+\rho (y-x)-z|^2}{4t}} \,\mathrm{d}\rho .
\end{align}
The Fubini theorem and (\ref{05132049.1}) leads to
\begin{align}
      & \int_{\mathbb{R}}|p_t(x,z)-p_t(y,z)|e^{\eta |z|}\,\mathrm{d}z \nonumber\\
      \leq &  \sqrt{\frac{2}{t}}|x-y| \int_0^1 \int_{\mathbb{R}} \frac{1}{\sqrt{2\pi\times 2t}} e^{-\frac{|x+\rho (y-x)-z|^2}{4t}} e^{\eta|z|} \,\mathrm{d}z \mathrm{d}\rho \nonumber\\
      \leq &  \sqrt{\frac{2}{t}}|x-y| \times \int_0^1 2 e^{\eta^2 t}\times e^{\eta |x+\rho (y-x)|} \,\mathrm{d}\rho  \nonumber\\
      \leq & 2\sqrt{2}\times \frac{|x-y|}{\sqrt{t}}\times e^{\eta^2 t}\times e^{\eta (|x|+|x-y|)} .
\end{align}

Proof of (iv). By (\ref{05161703.1}), the Fubini theorem and (\ref{05161718}), we have
{\allowdisplaybreaks \begin{align}
    & \int_{\mathbb{R}}|p_t(x,z)-p_t(y,z)|e^{\eta |z|}\eta |z| \,\mathrm{d}z \nonumber\\
    \leq & \sqrt{\frac{2}{t}}|x-y| \int_0^1 \int_{\mathbb{R}} \frac{1}{\sqrt{2\pi\times 2t}} e^{-\frac{|x+\rho (y-x)-z|^2}{4t}} e^{\eta|z|}\eta |z| \,\mathrm{d}z \mathrm{d}\rho \nonumber\\
    \leq & \sqrt{\frac{2}{t}}|x-y| \int_0^1 \bigg[e^{\eta^2 t}\times e^{\eta |x+\rho (y-x)|}\eta|x+\rho(y-x)| \nonumber\\
    & + 2e^{\eta^2 t}\bigg(2\eta^2 t+ \eta\sqrt{\frac{t}{\pi}}\bigg) e^{\eta |x+\rho (y-x)|}\bigg] \,\mathrm{d}\rho \nonumber\\
    \leq & \frac{\sqrt{2}|x-y|}{\sqrt{t}} \times \bigg[ e^{\eta^2 t}\times e^{\eta (|x|+|x-y|)}\eta (|x|+|x-y|) \nonumber\\
    & + 2e^{ \eta^2 t }\bigg(2\eta^2 t+ \eta\sqrt{\frac{t}{\pi}}\bigg)e^{\eta (|x|+|x-y|)}\bigg].
\end{align} }

Proof of (iv). This inequality can be found in Lemma 6.2 of \cite{S}, here we just give the explicit constant by straightforward calculations.

This completes the proof of Lemma \ref{A.4}.
$\blacksquare$

\section{Moment estimates}
\setcounter{equation}{0}
In this section, we will establish estimates for moments of stochastic convolutions, and obtain some a priori estimates  for solutions of equation (\ref{1.a}).

We begin with the  estimates of high order moments of stochastic convolutions.We stress that the precise lower order moment estimates are harder to get.
\begin{lemma}\label{210121.1036}
Let $h: \mathbb{R}_{+} \longmapsto \mathbb{R}_{+}$ be an increasing function.
Let $\{\sigma(s,y): (s,y)\in\mathbb{R}_+\times [0,1]\}$ be a random field such that the following stochastic convolution with respect to space time white noise is well defined.
Let $\tau$ be a stopping time.
%Assume that the following stochastic convolution with respect to stochastic field $\{\sigma(t,x): t\geq 0,x\in\mathbb{R}\}$ is well-defined.
Then  for any $p>10$ and $T>0$, there exists a constant $C_{p, h(T), T}>0$ such that
\begin{align}\label{210119.2121}
	& \mathbb{E} \sup_{(t,x)\in [0,T\wedge\tau]\times \mathbb{R}}\Bigg\{\left|\int_0^t\int_{\mathbb{R}}p_{t-s}(x,y)\sigma(s,y)\,W(\mathrm{d}s,\mathrm{d}y)\right| e^{-h(t)|x|}\Bigg\}^p \nonumber\\
	\leq & C_{p, h(T), T}\,{\mathbb{E}}\int_0^{T\wedge\tau} \int_{\mathbb{R}}|\sigma(t,x)|^p e^{-p h(t) |x|} \, \mathrm{d}x \mathrm{d}t .
\end{align}
%Then
%	\begin{align}\label{sto_con1}
%		\mathbb{E} \sup_{t\leq T, x\in\mathbb{R}}\Bigg\{\left|\int_0^T\int_{\mathbb{R}}p_{t-s}(x,y)\sigma(s,y)W(ds,dy)\right|^p e^{-\lambda|x|}\Bigg\} <\infty .
%	\end{align}
In particular, if $\sigma$ is bounded and $h$ is a positive constant, then the left hand side of (\ref{210119.2121}) is finite.
\end{lemma}
%\begin{remark}
%	The constant $C_{p, h(T), T}$ is given by (\ref{C_{T,p}}) and bounded by (\ref{210217.2153}).
%\end{remark}
\noindent {\bf Proof}.
%To prove the above lemma, it is sufficient to prove (\ref{sto_con1}) for $p>10$.
We employ the factorization method (see e.g. \cite{DPZ}). The proof here  is inspired by \cite{SZ1}.  Choose $\frac{3}{2p}<\alpha<\frac{1}{4}-\frac{1}{p}$. This is possible because we assume $p>10$. Let
\begin{align}
  (J_{\alpha}\sigma)(s,y):&= \int_0^s\int_{\mathbb{R}} (s-r)^{-\alpha}p_{s-r}(y,z)\sigma(r,z)\,W(\mathrm{d}r,\mathrm{d}z), \\
  (J^{\alpha-1}f)(t,x):&= \frac{\sin\pi\alpha}{\pi}\int_0^t\int_{\mathbb{R}} (t-s)^{\alpha-1}p_{t-s}(x,y)f(s,y)\,\mathrm{d}s\mathrm{d}y.
\end{align}

\noindent From the stochastic Fubini theorem (see Theorem 2.6 in \cite{W}),  it follows that for any $(t,x)\in\mathbb{R}_+\times \mathbb{R}$,
\begin{align}\label{103.1}
  \int_0^t\int_{\mathbb{R}} p_{t-s}(x,y)\sigma(s,y)\,W(\mathrm{d}s,\mathrm{d}y)=J^{\alpha-1}(J_{\alpha}\sigma)(t,x).
\end{align}
%In the sequel, we will use the following bounds, which can be seen by straightforward calculations.
%\begin{align}
%\label{05132049.1} \int_{\mathbb{R}} p_{t}(x,y)e^{\eta |y|} dy \leq & 2e^{\frac{\eta^2 t}{2}}e^{\eta |x|}, \quad \forall\ x,\eta\in\mathbb{R} , \\
%  \label{2102121.2100} \int_{\mathbb{R}} p_{t}(x,y)^2 e^{\eta |y|} dy \leq & \frac{1}{\sqrt{\pi t}}e^{\frac{\eta^2 t}{4}}e^{\eta |x|}, \quad \forall\ x,\eta\in\mathbb{R} .
%\end{align}

\noindent By H\"{o}lder's inequality and (\ref{2102121.2100}), we have
{\allowdisplaybreaks\begin{align}\label{104.1}
  & {\mathbb{E}}\sup_{(t,x)\in[0,T\wedge\tau]\times\mathbb{R}}\bigg\{\left|\int_0^t\int_{\mathbb{R}} p_{t-s}(x,y)\sigma(s,y)\,W(\mathrm{d}s,\mathrm{d}y)\right| e^{-h(t) |x|} \bigg\}^p \nonumber\\
  =& {\mathbb{E}}\sup_{(t,x)\in[0,T\wedge\tau]\times\mathbb{R}}\bigg\{\left|\frac{\sin\pi\alpha}{\pi} \int_0^t\int_{\mathbb{R}} (t-s)^{\alpha-1} p_{t-s}(x,y)J_{\alpha}\sigma(s,y)\,\mathrm{d}s\mathrm{d}y\right| e^{-h(t) |x|}\bigg\}^p \nonumber\\
  \leq & \left|\frac{\sin\pi\alpha}{\pi}\right|^p {\mathbb{E}}\sup_{(t,x)\in[0,T\wedge\tau]\times\mathbb{R}}\bigg\{\bigg[\int_0^t (t-s)^{\alpha-1} \nonumber\\
  &~~~~~~~~~~~~~\times\left(\int_{\mathbb{R}} p_{t-s}(x,y)|J_{\alpha}\sigma(s,y)|\,\mathrm{d}y\right)\,\mathrm{d}s\bigg]^p e^{-ph(t) |x|}\bigg\} \nonumber\\
  \leq & \left|\frac{\sin\pi\alpha}{\pi}\right|^p {\mathbb{E}}\sup_{(t,x)\in[0,T\wedge\tau]\times \mathbb{R}}\Bigg\{\bigg[\int_0^t (t-s)^{\alpha-1}  \nonumber\\
  &~~~~~~~~~~~~~\times\left(\int_{\mathbb{R}} p_{t-s}(x,y)e^{\frac{ph(t)}{2}|y|}e^{-\frac{ph(t)}{2}|y|}|J_{\alpha}\sigma(s,y)|^{\frac{p}{2}}\,\mathrm{d}y\right)^{\frac{2}{p}}\,\mathrm{d}s\bigg]^p e^{-p h(t) |x|}\Bigg\}\nonumber\\
  \leq & \left|\frac{\sin\pi\alpha}{\pi}\right|^p {\mathbb{E}}\sup_{(t,x)\in[0,T\wedge\tau]\times \mathbb{R}}\Bigg\{\bigg[\int_0^t (t-s)^{\alpha-1}  \left(\int_{\mathbb{R}} p_{t-s}(x,y)^2 e^{ph(t) |y|}\,\mathrm{d}y\right)^{\frac{1}{2}\times\frac{2}{p}} \nonumber\\
  &~~~~~~~~~~~~~\times\left(\int_{\mathbb{R}}|J_{\alpha}\sigma(s,y)|^p e^{-ph(t) |y|}\,\mathrm{d}y\right)^{\frac{1}{2}\times\frac{2}{p}}\,\mathrm{d}s\bigg]^p e^{-ph(t) |x|}\Bigg\} \nonumber\\
  \leq & \left|\frac{\sin\pi\alpha}{\pi}\right|^p \frac{1}{\sqrt{\pi}} e^{\frac{p^2h(T)^2}{4}T} \times {\mathbb{E}}\sup_{t\in[0,T\wedge\tau]} \Bigg[\int_0^t (t-s)^{\alpha-1-\frac{1}{2p}} \nonumber\\
  & ~~~~~~~~~~~~~~\times \left(\int_{\mathbb{R}}|J_{\alpha}\sigma(s,y)|^p e^{-ph(t) |y|}\,\mathrm{d}y\right)^{\frac{1}{p}}\,\mathrm{d}s\Bigg]^p  \nonumber\\
  \leq & \left|\frac{\sin\pi\alpha}{\pi}\right|^p \frac{1}{\sqrt{\pi}} e^{\frac{p^2h(T)^2}{4}T} \times {\mathbb{E}}\sup_{t\in[0,T\wedge\tau]}\Bigg\{ \left(\int_0^t (t-s)^{(\alpha-1-\frac{1}{2p})\frac{p}{p-1}}\,\mathrm{d}s\right)^{\frac{p-1}{p}\times p}  \nonumber\\
  &~~~~~~~~~~~~~\times\left(\int_0^t\int_{\mathbb{R}}|J_{\alpha}\sigma(s,y)|^p e^{-ph(t) |y|}\,\mathrm{d}y\mathrm{d}s\right)^{\frac{1}{p}\times p} \Bigg\} \nonumber\\
  \leq & \left|\frac{\sin\pi\alpha}{\pi}\right|^p \frac{1}{\sqrt{\pi}} e^{\frac{p^2h(T)^2}{4}T} \times \left(\int_0^T s^{(\alpha-1-\frac{1}{2p})\frac{p}{p-1}}\,\mathrm{d}s\right)^{p-1} \nonumber\\
  & ~~~~~~~~~~~~~\times\int_0^T\int_{\mathbb{R}} {\mathbb{E}}\left[|J_{\alpha}\sigma(s,y)|^p \mathbbm{1}_{[0, \tau]}(s) \right] e^{-ph(s) |y|}\,\mathrm{d}y\mathrm{d}s \nonumber\\
  \leq & C_{p, h(T), T, \alpha}^{\prime} \int_0^T\int_{\mathbb{R}} {\mathbb{E}}\Big[ \left|\int_0^s\int_{\mathbb{R}} (s-r)^{-\alpha}p_{s-r}(y,z)\sigma(r,z)\,W(\mathrm{d}r,\mathrm{d}z)\right|^p \nonumber\\
  &~~~~~~~~~~~~~~~~~~~~~~~~~~~~ \times\mathbbm{1}_{[0, \tau]}(s) \Big] e^{-ph(s) |y|}\,\mathrm{d}y\mathrm{d}s ,
%  \leq & C_{\lambda, p, T}^{\prime}  {\mathbb{E}}\int_0^T\int_{\mathbb{R}} |J_{\alpha}\sigma(s,y)|^p e^{-\lambda |y|}\,\mathrm{d}y\mathrm{d}s  ,
\end{align}}

\noindent where we have used the condition $\alpha >\frac{3}{2p}$ in the last inequality, so that
\begin{align}\label{210217.1615}
	C_{p, h(T), T, \alpha}^{\prime} =& \left|\frac{\sin\pi\alpha}{\pi}\right|^p \frac{1}{\sqrt{\pi}} e^{\frac{p^2h(T)^2}{4}T} \times \left(\int_0^T s^{(\alpha-1-\frac{1}{2p})\frac{p}{p-1}}\,\mathrm{d}s\right)^{p-1} \nonumber\\
	=&  \left|\frac{\sin\pi\alpha}{\pi}\right|^p \frac{1}{\sqrt{\pi}} e^{\frac{p^2h(T)^2}{4}T} \times \left(\frac{p-1}{\alpha p-\frac{3}{2}}\right)^{p-1} T^{\alpha p -\frac{3}{2}} .
\end{align}
For any fixed $(s,y)\in[0,T]\times\mathbb{R}$, let
\begin{align}
	M_t:= \int_0^t\int_{\mathbb{R}} (s-r)^{-\alpha}p_{s-r}(y,z)\sigma(r,z) \mathbbm{1}_{[0, \tau]}(r) \,W(\mathrm{d}r,\mathrm{d}z), \quad t\in [0,s] .
\end{align}
Then it is easy to see that $\{M_t\}_{t\in [0,s]}$ is a martingale. Applying the Bukrholder-Davis-Gundy inequality (see Proposition 4.4 in \cite{K} and also \cite{W}), we have for $t\in [0,s]$,
\begin{align}\label{210217.1647}
	\mathbb{E}|M_t|^p \leq & (4p)^{\frac{p}{2}} \,\mathbb{E}\langle M\rangle_t^{\frac{p}{2}} \nonumber\\
=&(4p)^{\frac{p}{2}}\,\mathbb{E}\left(\int_0^{t\wedge\tau}\int_{\mathbb{R}} (s-r)^{-2\alpha}p_{s-r}(y,z)^2\sigma(r,z)^2\,\mathrm{d}r\mathrm{d}z\right)^{\frac{p}{2}}.
\end{align}
Hence by the local property of the stochastic integral (see Lemma A.1 in Appendix of \cite{SZ1}),  we get
\begin{align}\label{210217.1647.1}
	 & {\mathbb{E}} \left[\left|\int_0^s\int_{\mathbb{R}} (s-r)^{-\alpha}p_{s-r}(y,z)\sigma(r,z)\,W(\mathrm{d}r,\mathrm{d}z)\right|^p \mathbbm{1}_{[0, \tau]}(s) \right]\nonumber\\
	 = & {\mathbb{E}} \left[\left|\int_0^s\int_{\mathbb{R}} (s-r)^{-\alpha}p_{s-r}(y,z)\sigma(r,z)\mathbbm{1}_{[0, \tau]}(r) \,W(\mathrm{d}r,\mathrm{d}z)\right|^p \mathbbm{1}_{[0, \tau]}(s) \right]\nonumber\\
	 \leq & {\mathbb{E}} \left|\int_0^s\int_{\mathbb{R}} (s-r)^{-\alpha}p_{s-r}(y,z)\sigma(r,z)\mathbbm{1}_{[0, \tau]}(r) \,W(\mathrm{d}r,\mathrm{d}z)\right|^p \nonumber\\
	 \leq & (4p)^{\frac{p}{2}}\,\mathbb{E}\left(\int_0^{s\wedge\tau}\int_{\mathbb{R}} (s-r)^{-2\alpha}p_{s-r}(y,z)^2\sigma(r,z)^2\,\mathrm{d}r\mathrm{d}z\right)^{\frac{p}{2}} .
\end{align}
%By the BDG inequality and H\"{o}ler's inequality, we have
Note that $p_t(x,y)\leq(2\pi t)^{-\frac{1}{2}}$ for any $t>0$ and $x,y\in\mathbb{R}$.  Using (\ref{210217.1647.1}), H\"{o}lder's inequality and the Fubini theorem, we get
\begin{align}\label{210217.1612}
% & {\mathbb{E}}\int_0^T\int_{\mathbb{R}} |J_{\alpha}\sigma(s,y)|^p e^{-\lambda |y|}\,\mathrm{d}y\mathrm{d}s \nonumber\\
%	=
	& \int_0^T\int_{\mathbb{R}} {\mathbb{E}}\left[ \left|\int_0^s\int_{\mathbb{R}} (s-r)^{-\alpha}p_{s-r}(y,z)\sigma(r,z)\,W(\mathrm{d}r,\mathrm{d}z)\right|^p \mathbbm{1}_{[0, \tau]}(s) \right] e^{-ph(s) |y|}\,\mathrm{d}y\mathrm{d}s \nonumber\\
	\leq & (4p)^{\frac{p}{2}}\int_0^T\int_{\mathbb{R}}{\mathbb{E}}\left\{\int_0^{s\wedge\tau}\int_{\mathbb{R}} (s-r)^{-2\alpha}p_{s-r}^2(y,z)\sigma^2(r,z)\, \mathrm{d}r \mathrm{d}z \right\}^{\frac{p}{2}} e^{-ph(s) |y|}\,\mathrm{d}y\mathrm{d}s \nonumber\\
	\leq & \left(\frac{4p}{\sqrt{2\pi}}\right)^{\frac{p}{2}}\int_0^T\int_{\mathbb{R}}{\mathbb{E}}\Bigg\{\left[\int_0^{s\wedge\tau}\int_{\mathbb{R}} \bigg|(s-r)^{-2\alpha -\frac{1}{2}}p_{s-r}^{1-\frac{2}{p}}(y,z)\bigg|^{\frac{p}{p-2}}\, \mathrm{d}r \mathrm{d}z\right]^{\frac{p-2}{p}} \nonumber\\
	&~~~~~~~~~~~\times \left[\int_0^{s\wedge\tau}\int_{\mathbb{R}} \bigg|p_{s-r}^{\frac{2}{p}}(y,z)\sigma^2(r,z)\bigg|^{\frac{p}{2}}\, \mathrm{d}r \mathrm{d}z\right]^{\frac{2}{p}}
	\Bigg\}^{\frac{p}{2}} e^{-ph(s) |y|}\,\mathrm{d}y\mathrm{d}s \nonumber\\
	\leq & \left(\frac{4p}{\sqrt{2\pi}}\right)^{\frac{p}{2}} \int_0^T  \left[\int_0^s  (s-r)^{-(2\alpha +\frac{1}{2}) \frac{p}{p-2} }\, \mathrm{d}r \right]^{\frac{p-2}{2}} \nonumber\\
	&~~~~~~~~~~~\times {\mathbb{E}} \int_0^{s\wedge\tau}\int_{\mathbb{R}} \left(\int_{\mathbb{R}}p_{s-r}(y,z)e^{-ph(s) |y|}\,\mathrm{d}y \right)  |\sigma(r,z)|^p \,\mathrm{d}z \mathrm{d}r\mathrm{d}s \nonumber\\
	\leq & C_{p, h(T), T, \alpha}^{\prime\prime} \,{\mathbb{E}} \int_0^{T\wedge\tau}\int_{\mathbb{R}} |\sigma(r,z)|^p e^{-ph(r)|z|} \, \mathrm{d}z \mathrm{d}r ,
\end{align}
where (\ref{05132049.1}) and the condition $\alpha<\frac{1}{4}-\frac{1}{p}$ are used to see that
\begin{align}\label{210217.1618}
	C_{p, h(T), T, \alpha}^{\prime\prime} = & \left(\frac{4p}{\sqrt{2\pi}}\right)^{\frac{p}{2}} \times T \times \left(\int_0^T  r^{-(2\alpha +\frac{1}{2}) \frac{p}{p-2} }\, \mathrm{d}r \right)^{\frac{p-2}{2}}  \times 2e^{\frac{p^2h(T)^2}{2}T}\nonumber\\
	= & \left(\frac{4p}{\sqrt{2\pi}}\right)^{\frac{p}{2}} \times \left(\frac{p-2}{\frac{p}{2}-2-2\alpha p}\right)^{\frac{p-2}{2}} T^{\frac{p}{4}-\alpha p}\times 2e^{\frac{p^2h(T)^2}{2}T}  .
\end{align}
Combining (\ref{104.1}) with (\ref{210217.1612}), we obtain
\begin{align}
  & \mathbb{E} \sup_{(t,x)\in[0,T\wedge\tau]\times \mathbb{R}}\Bigg\{\left|\int_0^t\int_{\mathbb{R}}p_{t-s}(x,y)\sigma(s,y)W(\mathrm{d}s,\mathrm{d}y)\right| e^{-h(t)|x|}\Bigg\}^p \nonumber\\
  \leq & C_{p, h(T), T} {\mathbb{E}} \int_0^{T\wedge\tau}\int_{\mathbb{R}} |\sigma(r,z)|^p e^{-ph(r)|z|} \, \mathrm{d}z \mathrm{d}r  ,
\end{align}
where
\begin{align}\label{C_{T,p}}
  C_{p, h(T), T}= \min_{\frac{3}{2p}<\alpha<\frac{1}{4}-\frac{1}{p}}C^{\prime}_{p, h(T), T, \alpha}\times C^{\prime\prime}_{p, h(T), T, \alpha} .
\end{align}
In view of (\ref{210217.1615}) and (\ref{210217.1618}), a straightforward calculation leads to
\begin{align}\label{210217.2153}
  C_{p, h(T), T} < 2 \sqrt{2}\,  p^{\frac{p}{2}}  \left(\frac{2}{\pi}\right)^p \left(\frac{1}{\sqrt{2\pi}}\right)^{\frac{p}{2}+1} \left(\frac{6p-8}{p-10}\right)^{\frac{3p}{2}-2}  T^{\frac{p}{4}-\frac{3}{2}} e^{\frac{3 p^2h(T)^2}{4}T}.
\end{align}
$\blacksquare$

\begin{proposition}\label{estimates 003}
Let $h: \mathbb{R}_{+} \longmapsto \mathbb{R}_{+}$ be an increasing function.
Let $\{\sigma(s,y): (s,y)\in\mathbb{R}_+\times [0,1]\}$ be a random field such that the following stochastic convolution with respect to the space time white noise is well defined.
%Let $\rho: \mathbb{R}_{+} \longmapsto \mathbb{R}_{+}$ be an increasing function. Assume that $\{\sigma(s,y): (s,y)\in\mathbb{R}_+\times [0,1]\}$ be a random field such that the stochastic integral against space time white noise is well defined.
Let $\tau$ be a stopping time.
Then for any $\epsilon, T>0$, and $0<p\leq 10$, there exists a constant $C_{\epsilon, p, h(T), T}$ such that
  \begin{align}\label{101.2}
    & \mathbb{E} \sup_{(t,x)\in [0,T\wedge\tau]\times \mathbb{R}}\Bigg\{\left|\int_0^t\int_{\mathbb{R}}p_{t-s}(x,y)\sigma(s,y)\,W(\mathrm{d}s,\mathrm{d}y)\right|  e^{-h(t)|x|}\Bigg\}^p \nonumber\\
    \leq & \epsilon \,\mathbb{E} \sup_{(t,x)\in[0,T\wedge\tau]\times\mathbb{R}}\left(|\sigma(t,x)| e^{-h(t)|x|} \right)^p  \nonumber\\
    & + C_{\epsilon,p,h(T),T} \,\mathbb{E}\int_0^{T\wedge\tau}\int_{\mathbb{R}} \left|\sigma(t,x)\right|^p e^{-ph(t)|x|}\,\mathrm{d}x\mathrm{d}t .
  \end{align}
\end{proposition}
\begin{remark}
	The constant $C_{\epsilon,p,h(T),T}$ is increasing with respect to $T$ and $C_{\epsilon,p,h(0),0}=0$ .
\end{remark}
\noindent {\bf Proof}. The following proof is  inspired by \cite{SZ1}. The proof is divided into two steps.

Step 1. We first show that for any $\rho, T>0$ and $q>10$,
  \begin{align}\label{102.1}
    & \mathbb{P}\left(\sup_{(t,x)\in[0,T\wedge\tau]\times\mathbb{R}}\left[\left|\int_0^t\int_{\mathbb{R}} p_{t-s}(x,y)\sigma(s,y)\,W(\mathrm{d}s,\mathrm{d}y) \right|e^{-h(t)|x|}\right]>\rho\right) \nonumber\\
    \leq & \mathbb{P}\left(\int_0^{T\wedge\tau}\int_{\mathbb{R}}\left|\sigma(s,y)\right|^q e^{-qh(s)|y|}\,\mathrm{d}y\mathrm{d}s >\rho^q\right) \nonumber\\
     & +  \frac{C_{q,h(T),T}}{\rho^q}\mathbb{E}\min\left\{\rho^q, \int_0^{T\wedge\tau}\int_{\mathbb{R}}\left|\sigma(s,y)\right|^q e^{-qh(t)|y|}\,\mathrm{d}y\mathrm{d}s\right\}.
  \end{align}
  Here the constant $C_{q,h(T), T}$ is the constant $C_{p,h(T),T}$ in (\ref{210119.2121}) with $p$ replaced by $q$.

To prove (\ref{102.1}), we set
\begin{align}
  \Omega_{\rho}:=\left\{\omega\in\Omega: \int_0^{T\wedge\tau} \int_{\mathbb{R}}|\sigma(s,y)|^q e^{-qh(s)|y|}\,\mathrm{d}y\mathrm{d}s\leq \rho^q\right\} .
\end{align}
By Chebyshev's inequality, we have
\begin{align}\label{105.1}
    & \mathbb{P}\left(\sup_{(t,x)\in[0,T\wedge\tau]\times\mathbb{R}}\left[\left|\int_0^t\int_{\mathbb{R}} p_{t-s}(x,y)\sigma(s,y)\,W(\mathrm{d}s,\mathrm{d}y) \right|e^{-h(t)|x|}\right]>\rho\right) \nonumber\\
    \leq & \mathbb{P}(\Omega\backslash\Omega_{\rho}) + \mathbb{P}\left(\sup_{(t,x)\in[0,T\wedge\tau]\times\mathbb{R}}\left[\left|\int_0^t\int_{\mathbb{R}} p_{t-s}(x,y)\sigma(s,y)\,W(\mathrm{d}s,\mathrm{d}y) \right| \mathbbm{1}_{\Omega_{\rho}} e^{-h(t)|x|} \right] >\rho\right) \nonumber\\
    \leq & \mathbb{P}\left(\int_0^{T\wedge\tau}\int_{\mathbb{R}}\left|\sigma(s,y)\right|^q e^{-qh(s)|y|}\,\mathrm{d}y\mathrm{d}s >\rho^q\right) \nonumber\\
    & + \frac{1}{\rho^q} \mathbb{E}\sup_{(t,x)\in[0,T\wedge\tau]\times\mathbb{R}} \left\{\left|\mathbbm{1}_{\Omega_{\rho}}\int_0^t\int_{\mathbb{R}} p_{t-s}(x,y)\sigma(s,y)\,W(\mathrm{d}s,\mathrm{d}y) \right|^q e^{-qh(t)|x|}\right\} .
\end{align}
%where, we have used Chebyshev's inequality in the last inequality, .
Now, we introduce the random field
\begin{align}
  \widetilde{\sigma}(s,y):= \sigma(s,y)\mathbbm{1}_{\left\{\omega\in\Omega: \ \int_0^{s\wedge\tau}\int_{\mathbb{R}}|\sigma(r,y)|^q e^{-qh(r)|y|} \,\mathrm{d}y\mathrm{d}r \leq\rho^q\right\}}.
\end{align}
Note that the stochastic integral of $\widetilde{\sigma}(\cdot, \cdot)$ with respect to the space time white noise $W$ is well defined.
%Since for any $\omega\in\Omega_{\beta}$,
%\begin{align}
%  \int_0^t\int_{\mathbb{R}} p_{t-s}(x,y)^2 |\sigma(s,y)-\widetilde{\sigma}(s,y)|^2\,\mathrm{d}s\mathrm{d}y =0 ,  \quad\forall\,t\in [0,T] ,
%\end{align}
By the local property of the stochastic integral (see Lemma A.1 in Appendix of \cite{SZ1}),
\begin{align}
  & \mathbbm{1}_{\Omega_{\rho}}\int_0^t\int_{\mathbb{R}} p_{t-s}(x,y)\sigma(s,y)\,W(\mathrm{d}s,\mathrm{d}y) \nonumber \\
  = & \mathbbm{1}_{\Omega_{\rho}}\int_0^t\int_{\mathbb{R}} p_{t-s}(x,y)\widetilde{\sigma}(s,y)\,W(\mathrm{d}s,\mathrm{d}y), \quad \mathbb{P}-a.s..
\end{align}
Hence using the bound (\ref{210119.2121}), we get
%{\allowdisplaybreaks
\begin{align}\label{105.2}
  & \mathbb{E}\sup_{(t,x)\in[0,T\wedge\tau]\times\mathbb{R}} \left\{\left|\mathbbm{1}_{\Omega_{\rho}}\int_0^t\int_{\mathbb{R}} p_{t-s}(x,y)\sigma(s,y)\,W(\mathrm{d}s,\mathrm{d}y) \right|^q e^{-qh(t)|x|}\right\} \nonumber\\
%  = & \mathbb{E}\sup_{(t,x)\in[0,T]\times\mathbb{R}} \left\{\left|\mathbbm{1}_{\Omega_{\beta}}\int_0^t\int_{\mathbb{R}} p_{t-s}(x,y)\widetilde{\sigma}(s,y)\,W(\mathrm{d}s,\mathrm{d}y) \right|^q e^{-q\lambda|x|}\right\} \nonumber\\
  \leq & \mathbb{E}\sup_{(t,x)\in[0,T\wedge\tau]\times\mathbb{R}} \left\{\left|\int_0^t\int_{\mathbb{R}} p_{t-s}(x,y)\widetilde{\sigma}(s,y)\,W(\mathrm{d}s,\mathrm{d}y) \right|^q e^{-qh(t)|x|}\right\} \nonumber\\
  \leq & C_{q,h(T),T} \,\mathbb{E}\int_0^{T\wedge\tau} \int_{\mathbb{R}} \left|\widetilde{\sigma}(s,y)\right|^q e^{-qh(s)|y|}\,\mathrm{d}s \nonumber\\
  \leq & C_{q,h(T),T} \,\mathbb{E}\min\left\{\rho^q, \int_0^{T\wedge\tau}\int_{\mathbb{R}}\left|\sigma(s,y)\right|^q e^{-qh(s)|y|}\,\mathrm{d}y\mathrm{d}s\right\} .
\end{align}
%}
Combining (\ref{105.1}) with (\ref{105.2}), we obtain (\ref{102.1}).

\vskip 0.3cm

Step 2.
%We show that for any $0<p\leq 10$, $q>10$ and $T>0$, there exists a constant $C_{p,q,h(t),T}$ such that
%  \begin{align}\label{102.2}
%    & \mathbb{E} \sup_{(t,x)\in[0,T]\times\mathbb{R}}\left\{\left|\int_0^t\int_{\mathbb{R}} p_{t-s}(x,y)\sigma(s,y)\,W(\mathrm{d}s,\mathrm{d}y) \right| e^{-\lambda|x|}  \right\}^p \nonumber\\
%    \leq & C_{p,q,h(t),T} \mathbb{E}\left[\int_0^T\int_{\mathbb{R}}\left|\sigma(s,y)\right|^q e^{-q\lambda|y|}\,\mathrm{d}y\mathrm{d}s \right]^{\frac{p}{q}} .
%  \end{align}
Let now $0<p\leq 10$. From (\ref{102.1}) and Lemma A.2 in Appendix of \cite{SZ1}, it follows that
%The estimate (\ref{102.2}) can be easily derived from (\ref{102.1}) and Lemma ??? in Appendix as follows:
{\allowdisplaybreaks
\begin{align}\label{210217.2157}
  & \mathbb{E} \sup_{(t,x)\in[0,T\wedge\tau]\times\mathbb{R}}\left\{\left|\int_0^t\int_{\mathbb{R}} p_{t-s}(x,y)\sigma(s,y)\,W(\mathrm{d}s,\mathrm{d}y) \right| e^{-h(t)|x|}  \right\}^p \nonumber\\
  = & \int_0^{\infty} p\rho^{p-1} \mathbb{P}\left(\sup_{(t,x)\in[0,T\wedge\tau]\times\mathbb{R}}\left[\left|\int_0^t\int_{\mathbb{R}} p_{t-s}(x,y)\sigma(s,y)\,W(\mathrm{d}s,\mathrm{d}y) \right| e^{-h(t)|x|}\right]>\rho\right)\,\mathrm{d}\rho \nonumber\\
  \leq & \int_0^{\infty} p\rho^{p-1} \mathbb{P}\left(\int_0^{T\wedge\tau} \int_{\mathbb{R}}|\sigma(s,y)|^q e^{-qh(s)|y|}\,\mathrm{d}y\mathrm{d}s>\rho^q\right)\,\mathrm{d}\rho \nonumber\\
  & + C_{q,h(T),T}\int_0^{\infty} p\rho^{p-1-q}\,\mathbb{E}\min\left\{\rho^q, \int_0^{T\wedge\tau} \int_{\mathbb{R}}|\sigma(s,y)|^q e^{-qh(s)|y|}\,\mathrm{d}y\mathrm{d}s\right\}\,\mathrm{d}\rho \nonumber\\
 = & C_{p,q,h(T),T}\,\mathbb{E} \left[\int_0^{T\wedge\tau}\int_{\mathbb{R}} |\sigma(s,y)|^q e^{-qh(s)|y|}\,\mathrm{d}y\mathrm{d}s\right]^{\frac{p}{q}} \nonumber\\
 \leq & C_{p,q,h(T),T} \,\mathbb{E}\Bigg[\sup_{(s,y)\in[0,T\wedge\tau]\times\mathbb{R}}\left(|\sigma(s,y)|e^{-h(s)|y|}\right)^{\frac{(q-p)p}{q}} \nonumber\\
  & ~~~~~~~~~~~~~~~~\times \left(\int_0^{T\wedge\tau}\int_{\mathbb{R}}\left|\sigma(s,y)\right|^p e^{-ph(s)|y|}\,\mathrm{d}y \mathrm{d}s \right)^{\frac{p}{q}} \Bigg] \nonumber\\
%    = & C_{T,p,q} \mathbb{E}\left[\sup_{(s,y)\in[0,T]\times[0,1]}|\sigma(s,y)|^{\frac{(q-p)p}{q}} \times\left(\int_0^T\sup_{y\in[0,1]}\left|\sigma(s,y)\right|^p\,ds\right)^{\frac{p}{q}} \right] \nonumber\\
    \leq & \epsilon \,\mathbb{E} \sup_{(s,y)\in[0,T\wedge\tau]\times \mathbb{R}}\left(|\sigma(s,y)|e^{-h(s)|y|}\right)^p  \nonumber\\
    & + C_{p,q,h(T),T}\times C_{\epsilon,p,q,h(T),T}\, \mathbb{E}\int_0^{T\wedge\tau}\int_{\mathbb{R}}\left|\sigma(s,y)\right|^p e^{-ph(s)|y|}\,\mathrm{d}y \mathrm{d}s ,
\end{align}}

\noindent where
\begin{align}\label{C_Tpq}
  C_{p,q,h(T),T}:=1+ \frac{q}{q-p} C_{q,h(T),T}  ,
\end{align}
and we have used the following Young inequality
\begin{align}\label{C_Tpqe}
  ab\leq & \frac{\epsilon}{C_{p,q,h(T),T}}\,a^{\frac{q}{q-p}} + C_{\epsilon,p,q,h(T),T}\,b^{\frac{q}{p}}, \quad \forall\, a,b>0, \nonumber\\  C_{\epsilon,p,q,h(T),T}:= & p\left(\frac{q-p}{\epsilon/C_{p,q,h(T),T}}\right)^{\frac{q-p}{p}} q^{-\frac{q}{p}} .
\end{align}
%and take
% \begin{align}
%   C_{\epsilon,p,q,T}= p\left(\frac{(q-p)C_{T,p,q}}{\epsilon}\right)^{\frac{q-p}{p}} q^{-\frac{q}{p}} ,
%%   \quad \varepsilon=\frac{\epsilon}{C_{T,p,q}} .
% \end{align}
Set
\begin{align}\label{C_{T,p,epsilon}}
  C_{\epsilon,p,h(T),T}:=\inf_{q>10} C_{p,q,h(T),T}\times C_{\epsilon,p,q,h(T),T} .
\end{align}

\noindent Combining (\ref{C_Tpq}) and (\ref{C_Tpqe}) gives
\begin{align}\label{C_Tpe}
	C_{\epsilon,p,h(T),T}= \inf_{q>10} \left\{\frac{p}{q-p} q^{-\frac{q}{p}} \epsilon^{1-\frac{q}{p}} \left(q-p+qC_{q,h(T),T}\right)^{\frac{q}{p}}\right\} ,
\end{align}
where the constant $C_{q,h(T),T}$ is bounded by the right hand side of (\ref{210217.2153}) with $p$ replaced by $q$.
Now, (\ref{101.2}) follows from (\ref{210217.2157}) with the constant $C_{\epsilon,p,h(T),T}$ defined above.
$\blacksquare$

\vskip 0.3cm

Next, we will establish an a priori estimate of solutions to (\ref{1.a}).  Throughout this paper, we will use the following  notations. For $\lambda, \kappa >0$, set
\begin{align}
	\beta(\lambda, \kappa):= & \max\left\{\frac{\lambda^2}{2}, 4\kappa \right\}, \\
	T^*(\lambda, \kappa):= & \frac{1}{2\beta(\lambda,\kappa)}\left[1+\log\left(\frac{4\beta(\lambda,\kappa)}{\lambda^2}\log\frac{\beta(\lambda,\kappa)}{2\kappa}\right)\right] .
\end{align}
It is easy to see that for any $\kappa>0$, $T^*(\lambda,\kappa)\rightarrow\infty$ as $\lambda\rightarrow 0$.

\begin{lemma}\label{lemma 3.1}
Assume that (H1) is satisfied and $\sigma$ is bounded.
Let $u$ be a solution of (\ref{1.a}).
%Let $\lambda$ be any positive number,
%$\beta:=\frac{\lambda^2}{2}\vee 4c_1$.
Set also
\begin{gather}
%\label{2.1}  T^*(\lambda):= \frac{1}{2\beta}\left[1+\log\left(\frac{4\beta}{\lambda^2}\log\frac{\beta}{2c_1}\right)\right] , \\
\label{2.2}  V(t,x):=\int_0^t\int_{\mathbb{R}} p_{t-s}(x,y)\sigma(u(s,y)) \,W(\mathrm{d}s,\mathrm{d}y) .
\end{gather}
Then for any $\lambda>0$ and $T\leq T^*(\lambda, c_1)$, there exists a constant $C_{\lambda,c_1,T}$ such that the following a priori estimate holds for $\mathbb{P}$-a.s.,
\begin{align}\label{2.3}
& \sup_{t\leq T, x\in\mathbb{R}}\left(|u(t,x)|e^{-\lambda |x|e^{\beta t}}\right) \nonumber\\
  \leq & C_{\lambda, c_1, T}\times \bigg\{1+ 2c_2 T +4e^{\frac{\lambda^2 T}{2}}\sup_{x\in\mathbb{R}}\left(|u_0(x)|e^{-\lambda |x|}\right) \nonumber\\
  & + 2   \sup_{(t,x)\in[0,T]\times\mathbb{R}}\left(|V(t,x)|e^{-\lambda |x|}\right)\bigg\}^{e^{4c_1 T e^{\frac{\lambda^2}{4\beta} e^{2\beta T -1}}}} ,
\end{align}
where we write $\beta$ instead of $\beta(\lambda, c_1)$ for simplicity, and the constant $c_1$ is same as that in condition (H1).
\end{lemma}
\begin{remark}
Lemma \ref{lemma 3.1} actually implies that the solutions of (\ref{1.a}) don't blow up in the space $C_{tem}$, since we can take sufficiently small $\lambda>0$ such that $T^*(\lambda,c_1)$ can be larger than any given number.
\end{remark}

\vskip 0.3cm
\noindent {\bf Proof}.
Set
\[U(T):= \sup_{(t,x)\in [0,T]\times\mathbb{R}}\left(|u(t,x)|e^{-\lambda |x|e^{\beta t}}\right) .\]
From (\ref{definition solution}), we have
\begin{align}\label{200404.1559}
	U(T) \leq & \sup_{t\leq T, x\in\mathbb{R}}\left(|P_t u_0(x)|e^{-\lambda |x|e^{\beta t}}\right) +  \sup_{t\leq T, x\in\mathbb{R}}\left(|V(t,x)|e^{-\lambda |x|}\right) \nonumber\\
	& + \sup_{t\leq T, x\in\mathbb{R}}\left\{\left|\int_0^t\int_{\mathbb{R}}p_{t-s}(x,y)b(u(s,y)) \,\mathrm{d}s\mathrm{d}y\right|e^{-\lambda |x|e^{\beta t}}\right\} .
\end{align}

\noindent Now we estimate the three terms on the right hand side of the above inequality.
\begin{align}\label{05161111}
  & \sup_{t\leq T, x\in\mathbb{R}}\left(|P_t u_0(x)|e^{-\lambda |x|e^{\beta t}}\right) \nonumber\\
  \leq & \sup_{t\leq T, x\in\mathbb{R}}\left\{\left|\int_{\mathbb{R}}p_t(x,y)u_0(y)\,\mathrm{d}y\right|e^{-\lambda |x|}\right\} \nonumber\\
  \leq & \sup_{y\in\mathbb{R}}\left(|u_0(y)|e^{-\lambda |y|}\right)\times \sup_{t\leq T, x\in\mathbb{R}}\left\{\int_{\mathbb{R}}p_t(x,y) e^{\lambda |y|}\,\mathrm{d}y \times e^{-\lambda |x|}\right\} \nonumber\\
  \leq & 2e^{\frac{\lambda^2 T}{2}} \sup_{y\in\mathbb{R}}\left(|u_0(y)|e^{-\lambda |y|}\right) ,
\end{align}
where we have used (\ref{05132049.1}).
%\begin{align}\label{05132049.1}
%  \int_{\mathbb{R}} p_{t}(x,y)e^{\eta |y|} dy \leq 2e^{\frac{\eta^2 t}{2}}e^{\eta |x|},
%\end{align}
%for any $x,\eta\in\mathbb{R}$.
Applying Lemma \ref{210121.1036} and using the boundedness of $\sigma$, we get that
for any $p, q, T>0$,
\begin{align}\label{05142145}
  \mathbb{E}\sup_{t\leq T, x\in\mathbb{R}}\left(|V(t,x)|^p e^{-q |x|}\right)<\infty .
\end{align}
In particular,
\begin{align}\label{200404.1600}
\sup_{t\leq T, x\in\mathbb{R}}\left(|V(t,x)|e^{-\lambda |x|}\right)<+\infty 	, \quad \mathbb{P}-a.s..
\end{align}

\noindent The nonlinear term can be estimated as follows. By (\ref{210124.2000}) and $\log_{+}(ab)\leq \log_{+}a +\log_{+}b$ for any $a,b>0$, we have
\begin{align}\label{05161114}
&
\sup_{t\leq T, x\in\mathbb{R}}\left\{\left|\int_0^t\int_{\mathbb{R}}p_{t-s}(x,y)b(u(s,y))\,\mathrm{d}s\mathrm{d}y\right|e^{-\lambda |x|e^{\beta t}}\right\} \nonumber\\
\leq &
\sup_{t\leq T, x\in\mathbb{R}}\left\{\int_0^t\int_{\mathbb{R}}p_{t-s}(x,y)\left( c_1|u(s,y)|\log_{+} |u(s,y)| + c_2\right)\,\mathrm{d}s\mathrm{d}y \times e^{-\lambda |x|e^{\beta t}}\right\} \nonumber\\
\leq &
c_2 T + c_1\sup_{t\leq T, x\in\mathbb{R}}\bigg\{\int_0^t\int_{\mathbb{R}}p_{t-s}(x,y) e^{\lambda |y| e^{\beta s}}\times\left(|u(s,y)|e^{-\lambda |y|e^{\beta s}}\right) \nonumber\\
&~~~~~~~~~~~~~~~~~\times\log_{+}\left[\left(|u(s,y)|e^{-\lambda |y| e^{\beta s}}\right)\times e^{\lambda |y| e^{\beta s}}\right] \,\mathrm{d}s\mathrm{d}y \times e^{-\lambda |x|e^{\beta t}}\bigg\} \nonumber\\
\leq &
c_2 T
+ c_1\sup_{t\leq T, x\in\mathbb{R}}\bigg\{\int_0^t \sup_{y\in\mathbb{R}} \left[\left(|u(s,y)|e^{-\lambda |y|e^{\beta s}} \right)\times \log_{+}\left(|u(s,y)|e^{-\lambda |y| e^{\beta s}}\right)\right] \nonumber\\
&~~~~~~~~~~~~~~~~\times \int_{\mathbb{R}}p_{t-s}(x,y) e^{\lambda |y| e^{\beta s}} \,\mathrm{d}y\mathrm{d}s \times e^{-\lambda |x|e^{\beta t}}\bigg\} \nonumber\\
&
+ c_1\sup_{t\leq T, x\in\mathbb{R}}\bigg\{\int_0^t \sup_{y\in\mathbb{R}}\left(|u(s,y)|e^{-\lambda |y|e^{\beta s}}\right) \nonumber\\
&~~~~~~~~~~~~~~~\times\int_{\mathbb{R}}p_{t-s}(x,y) e^{\lambda |y| e^{\beta s}} \lambda |y| e^{\beta s} \,\mathrm{d}y\mathrm{d}s \times e^{-\lambda |x|e^{\beta t}}\bigg\} \nonumber\\
=& :  c_2 T + I + II.
\end{align}
Note that the function $x\mapsto x\log_{+}x$ is increasing on $[0,\infty)$, so we have
\begin{align}\label{05161112}
I \leq &
c_1\sup_{t\leq T, x\in\mathbb{R}}\bigg\{\int_0^t \sup_{y\in\mathbb{R},r\leq s} \left[\left(|u(r,y)|e^{-\lambda |y|e^{\beta r}} \right)\times \log_{+}\left(|u(r,y)|e^{-\lambda |y| e^{\beta r}}\right)\right] \nonumber\\
& ~~~~~~~~~~~~~~\times 2 e^{\frac{\lambda^2 (t-s) e^{2\beta s}}{2}} e^{\lambda |x|e^{\beta s}} \,\mathrm{d}s \times e^{-\lambda |x|e^{\beta t}}\bigg\} \nonumber\\
\leq & 2c_1 \sup_{t\leq T}\bigg\{ \sup_{s\leq t}\left(e^{\frac{\lambda^2 (t-s) e^{2\beta s}}{2}}\right) \int_0^t U(s)\log_{+} U(s) \,\mathrm{d}s \bigg\} \nonumber\\
\leq & 2c_1  e^{\frac{\lambda^2}{4\beta}e^{2\beta T-1}} \int_0^T U(s)\log_{+} U(s) \,\mathrm{d}s ,
\end{align}

\noindent where we have used (\ref{05132049.1}) and
\begin{align}\label{05132049.2}
  \max_{s\in [0,t]} e^{\frac{\lambda^2 (t-s) e^{2\beta s}}{2}} =e^{\frac{\lambda^2}{4\beta}e^{2\beta t-1}}.
\end{align}
For the term $II$, we estimate as follows
{\allowdisplaybreaks \begin{align}\label{05132036}
  II \leq & c_1\sup_{t\leq T, x\in\mathbb{R}}\bigg\{\int_0^t \sup_{y\in\mathbb{R}}\left(|u(s,y)|e^{-\lambda |y|e^{\beta s}}\right) \nonumber\\
&~~~~~~~~~\times \left(e^{\frac{\lambda^2 (t-s) e^{2\beta s}}{2}}e^{\lambda |x| e^{\beta s}} \lambda |x| e^{\beta s} + C_{\lambda,\beta, t}e^{\lambda |x|e^{\beta s}} \right) \mathrm{d}s \times e^{-\lambda |x|e^{\beta t}}\bigg\} \nonumber\\
\leq & c_1\sup_{t\leq T, x\in\mathbb{R}}\bigg\{\sup_{s\leq t, y\in\mathbb{R}}\left(|u(s,y)|e^{-\lambda |y|e^{\beta s}}\right) \times \frac{1}{\beta}\sup_{s\leq t}\left(e^{\frac{\lambda^2 (t-s) e^{2\beta s}}{2}}\right)\nonumber\\
&~~~~~~~~~~\times \int_0^t \frac{\mathrm{d}}{\mathrm{d}s} e^{\lambda |x| e^{\beta s}} \,\mathrm{d}s\times e^{-\lambda |x|e^{\beta t}}\bigg\} \nonumber\\
& + c_1\sup_{t\leq T}\bigg\{C_{\lambda,\beta, t} \int_0^t \sup_{r\leq s, y\in\mathbb{R}}\left(|u(r,y)|e^{-\lambda |y|e^{\beta r}}\right) \mathrm{d}s \bigg\} \nonumber\\
\leq & \frac{c_1}{\beta}e^{\frac{\lambda^2}{4\beta}e^{2\beta T-1}} U(T) + C_{\lambda, c_1, T}\int_0^T U(s)\,\mathrm{d}s ,
\end{align}
}

\noindent where we have used (\ref{05161718}),
%estimate
%\begin{align}\label{05161718}
%\int_{\mathbb{R}}p_{t}(x,y) e^{\eta |y|} \eta |y|  dyds \leq e^{\frac{\eta^2 t}{2}}e^{\eta |x|}\eta|x| + 2e^{\frac{\eta^2 t}{2}}\left(\eta^2 t+ \eta\sqrt{\frac{t}{2\pi}}\right) e^{\eta |x|} ,
%\end{align}
(\ref{05132049.2}), and that the constant $C_{\lambda,\beta,t}$ is increasing with respect to  $t>0$. Note that
\begin{align}\label{05132040.1}
  \frac{c_1}{\beta}e^{\frac{\lambda^2}{4\beta}e^{2\beta T-1}}\leq \frac{1}{2} \iff T\leq  T^*(\lambda, c_1)= \frac{1}{2\beta}\left[1+\log\left(\frac{4\beta}{\lambda^2}\log\frac{\beta}{2c_1}\right)\right] .
\end{align}
Hence for $T\leq T^*(\lambda, c_1)$,
\begin{align}\label{05132040.2}
  II \leq & \frac{1}{2} U(T) +  c_1 C_{\lambda, \beta, T}\int_0^T U(s)\,\mathrm{d}s .
\end{align}

\noindent Combining (\ref{200404.1559}), (\ref{05161111}), (\ref{200404.1600}), (\ref{05161114}), (\ref{05161112}) and (\ref{05132040.2}) together, we obtain that for $T\leq T^*(\lambda, c_1)$,
\begin{align}
  U(T)\leq & 2e^{\frac{\lambda^2 T}{2}} \sup_{y\in\mathbb{R}}\left(|u_0(y)|e^{-\lambda |y|}\right) +  \sup_{t\leq T, x\in\mathbb{R}}\left(|V(t,x)| e^{-\lambda |x|}\right) \nonumber\\
  & + c_2 T +2c_1 e^{\frac{\lambda^2}{4\beta}e^{2\beta T-1}} \int_0^T U(s)\log_{+} U(s) \,\mathrm{d}s  \nonumber\\
  & + \frac{1}{2} U(T) +  C_{\lambda, c_1, T}\int_0^T U(s)\,\mathrm{d}s .
\end{align}
Subtracting $\frac{1}{2}U(T)$ on both sides of the above inequality, and then applying the log Gronwall inequality (see Lemma \ref{A.1}), (\ref{2.3}) is deduced.
$\blacksquare$

\section{Existence of weak solutions}
\setcounter{equation}{0}

In this section, we assume that (H1) is satisfied and that  $\sigma$ is bounded, continuous.
We will  approximate the coefficients $b$ and $\sigma$ by  Lipschitz continuous  functions and establish the existence of  weak solutions of the stochastic reaction-diffusion equation.

% hence deduce strong solutions solutions in the probabilistic sense by Yamada-Watanabe theorem.
%Next, we prove the existence of the solutions.
%By the Yamada-Watanabe theorem, it suffices to prove the existence of the weak solutions in the probabilistic sense.
%Set $b(x):= x\log|x|$.

Let $\varphi$ be a nonnegative smooth function on $\mathbb{R}$ such that the support of $\varphi$ is contained in $(-1,1)$ and $\int_{\mathbb{R}}\varphi(x)\,\mathrm{d}x=1$. Let $\{\eta_n\}_{n\geq 1}$ be a sequence of symmetric smooth functions such that for any $n\geq 1$, $0\leq \eta_n \leq 1$, $\eta_n(x)=1$ if $|x|\leq n$, and $\eta_n(x)=0$ if $|x|\geq n+2$.
Define
\begin{align}
  b_n(x):=n\int_{\mathbb{R}}b(y)\varphi(n(x-y))\,\mathrm{d}y \times \eta_n(x) , \\
  \sigma_n(x):=n\int_{\mathbb{R}}\sigma(y)\varphi(n(x-y))\,\mathrm{d}y \times \eta_n(x) .
\end{align}
Assume that $\sigma$ is bounded by a constant $K_{\sigma}$, that is
\begin{align}
	|\sigma(z)|\leq K_{\sigma}, \quad \forall\, z\in\mathbb{R}.
\end{align}
Then it is easy to check that there exist constants $L_n, L_b$ and $K_n$ such that for any $x,y\in\mathbb{R}$,
\begin{align}
\label{05160929.1}  |b_n(x)-b_n(y)|\leq &  L_n|x-y| , \\
\label{05160929.2}  |b_n(x)|\leq & c_1 |x|\log_{+}|x| + L_b(|x|+1 ), \\
\label{210122.1854}  |\sigma_n(x)-\sigma_n(y)|\leq &  K_n|x-y| , \\
\label{210122.1855}  |\sigma_n(x)|\leq & K_{\sigma} ,
\end{align}
where the constant $c_1$ is same as that in condition (H1), and the constant $L_b$ is independent of $n$. Moreover, if $x_n\rightarrow x$ in $\mathbb{R}$, then
\begin{align}
\label{05160930.1}
  b_n(x_n)\rightarrow b(x) , \\
  \label{210126.2140}
  \sigma_n(x_n)\rightarrow \sigma(x) .
\end{align}

For $n\geq 1$, consider the following stochastic equation on the real-line $\mathbb{R}$,
\begin{align}\label{05150908}
  u_n(t,x)= & P_t u_0(x)+\int_0^t\int_{\mathbb{R}} p_{t-s}(x,y) b_n(u_n(s,y)) \,\mathrm{d}s\mathrm{d}y \nonumber\\
   & + \int_0^t\int_{\mathbb{R}} p_{t-s}\sigma_n(u_n(s,y)) \,W(\mathrm{d}s,\mathrm{d}y).
\end{align}
It is known (see \cite{S,MPS,MP}) that for each $n\geq 1$, there exists a unique solution $u_n$ to the above equation. Moreover, the sample paths of $u_n$ are a.s. in $C(\mathbb{R}_{+},C_{tem})$.
The following result is a uniform bound for the solutions $u_n$.
\begin{lemma}\label{lemma 05151406}
Assume  $u_0 \in C_{tem}$ and (H1). Suppose that  $\sigma$ is bounded and continuous.
Then for any $p\geq 1$ and $\lambda, T >0$, we have
  \begin{align}\label{05142144}
    \sup_{n\geq 1} \mathbb{E}\left[\sup_{t\leq T, x\in\mathbb{R}} \left(|u_n(t,x)|e^{-\lambda |x|}\right)^p\right] <\infty .
  \end{align}
\end{lemma}
\noindent {\bf Proof}. It suffices to prove this lemma for sufficiently large $p$ and sufficiently small $\lambda$. Fix $T>0$.  As $\lim_{\lambda\rightarrow 0}T^*(\lambda, c_1)=\infty$, there exists a positive constant $\lambda_T$  such that $T\leq T^*(\lambda, c_1)$ for all $\lambda\leq \lambda_T$. For any fixed $\lambda\leq\lambda_T$, we choose $\lambda_0 >0$ so that $2\lambda_0 e^{\beta(\lambda_0, c_1) T} = \lambda$.
In the following we write $\beta$ for $\beta(\lambda_0, c_1)$ to simplify the notation.
Define
\begin{align}
  U_n(r):=\sup_{t\leq r, x\in\mathbb{R}}\left(|u_n(t,x)|e^{-\lambda_0 |x|e^{\beta t}}\right),\quad r\in[0, T].
\end{align}
It remains  to prove
\begin{align}\label{210225.1905}
	\sup_{n\geq 1}\mathbb{E} [U_n(T)]^p <\infty .
\end{align}
From (\ref{05150908}) we have
\begin{align}
  U_n(r) \leq & \sup_{t\leq r, x\in\mathbb{R}} \left(|P_t u_0(x)|e^{-\lambda_0 |x| e^{\beta t}}\right) \nonumber\\
  & + \sup_{t\leq r, x\in\mathbb{R}}\left\{\left|\int_0^t\int_{\mathbb{R}}p_{t-s}(x,y)b_n(u_n(s,y))\,\mathrm{d}s\mathrm{d}y\right|e^{-\lambda_0 |x|e^{\beta t}}\right\} \nonumber\\
  & + \sup_{t\leq r, x\in\mathbb{R}}\left\{\left|\int_0^t\int_{\mathbb{R}} p_{t-s}(x,y) \sigma_n(u_n(s,y))\,W(\mathrm{d}s,\mathrm{d}y)\right|e^{-\lambda_0 |x|e^{\beta t}}\right\} .
\end{align}
In the above inequality,
the first term can be estimated the same as (\ref{05161111}).
Let
\begin{align}
	V_n(t,x) = \int_0^t\int_{\mathbb{R}} p_{t-s}(x,y) \sigma_n(u_n(s,y))\,W(\mathrm{d}s,\mathrm{d}y) .
\end{align}
Then by (\ref{210122.1855}) and Lemma \ref{210121.1036}, we have for any $p>0$,
\begin{align}\label{210124.2114}
	 & \sup_{n\geq 1}\mathbb{E} \sup_{t\leq T, x\in\mathbb{R}}\left\{\left|V_n(t,x)\right|e^{-\lambda_0 |x|e^{\beta t}}\right\}^p  \nonumber\\
	\leq & \sup_{n\geq 1}\mathbb{E}\sup_{t\leq T, x\in\mathbb{R}}\left\{\left|V_n(t,x)\right|e^{-\lambda_0 |x| }\right\}^p \nonumber\\
	\leq & C_{\lambda_0, K_{\sigma}, T, p} <\infty.
\end{align}
%and for any $p>0$,
%\begin{align}
%	\sup_{n\geq 1}\mathbb{E} \bigg\{\sup_{t\leq T, x\in\mathbb{R}}\left|V_n(t,x) e^{-\lambda_0 |x|}\right|^p \bigg\} <\infty ,
%\end{align}
%where we have used (\ref{210122.1855}) and Lemma \ref{210121.1036}.
On the other hand,
\begin{align}\label{4.0}
  & \sup_{t\leq r, x\in\mathbb{R}}\left\{\left|\int_0^t\int_{\mathbb{R}}p_{t-s}(x,y)b_n(u_n(s,y))\,\mathrm{d}s\mathrm{d}y\right|e^{-\lambda_0 |x|e^{\beta t}}\right\} \nonumber\\
  \leq & \sup_{t\leq r, x\in\mathbb{R}}\bigg\{\int_0^t\int_{\mathbb{R}}p_{t-s}(x,y)[ c_1|u_n(s,y)|\log_+ |u_n(s,y)| \nonumber\\
  & ~~~~~~~~~~~~ +L_b(|u_n(s,y)|+1) ]\,\mathrm{d}s\mathrm{d}y\times e^{-\lambda_0 |x|e^{\beta t}}\bigg\} \nonumber\\
 \leq & c_1\sup_{t\leq r, x\in\mathbb{R}}\left\{\int_0^t\int_{\mathbb{R}}p_{t-s}(x,y)|u_n(s,y)|\log_+ |u_n(s,y)|\,\mathrm{d}s\mathrm{d}y\times e^{-\lambda_0 |x|e^{\beta t}}\right\} \nonumber\\
 & + L_b \sup_{t\leq r, x\in\mathbb{R}}\left\{\int_0^t\int_{\mathbb{R}}p_{t-s}(x,y)  |u_n(s,y)|\,\mathrm{d}s\mathrm{d}y\times e^{-\lambda_0 |x|e^{\beta t}}\right\} + L_b r .
% \nonumber\\
% & + Lr+ \sup_{t\leq r, x\in\mathbb{R}}\left\{\int_0^t\int_{\mathbb{R}}p_{t-s}(x,y)|u_n(s,y)|dsdy\times e^{-\lambda |x|e^{\beta t}}\right\} .
\end{align}

\noindent Now using (\ref{4.0}) and following a similar proof as that of Lemma \ref{lemma 3.1} we obtain
\begin{align}
  U_n(T)\leq & C_{\lambda_0, c_1, L_b, T}\times \bigg\{1+2L_b T+4e^{\frac{\lambda_0^2 T}{2}}\sup_{x\in\mathbb{R}}\left(|u_0(x)|e^{-\lambda_0 |x|}\right) \nonumber\\
  & + 2\sup_{t\leq T, x\in\mathbb{R}}\left(|V_n(t,x)|e^{-\lambda_0 |x|}\right)\bigg\}^{e^{4c_1 T e^{\frac{\lambda_0^2}{4\beta} e^{2\beta T -1}}}} .
\end{align}
%Note that the right hand side of the above inequality don't dependent on $n$.
Hence  it follows from (\ref{210124.2114}) that
\begin{align}\label{210225.1854}
  \mathbb{E} [U_n(T)^p] \leq & C_{\lambda_0, c_1, L_b, T, p}\left[1+\sup_{x\in\mathbb{R}}\left(|u_0(x)|e^{-\lambda_0 |x|}\right)^{p e^{4c_1 T e^{\frac{\lambda_0^2}{4\beta} e^{2\beta T -1}}}}\right] \nonumber\\\
  & +  C_{\lambda_0, c_1, L_b, T, p} \,\mathbb{E}\left[\sup_{t\leq T, x\in\mathbb{R}}\left(\left|V_n(t,x)\right|e^{-\lambda_0 |x| }\right)^{p e^{4c_1 T e^{\frac{\lambda_0^2}{4\beta} e^{2\beta T -1}}}}\right] \nonumber\\
  \leq &  C_{\lambda_0, c_1, L_b, K_{\sigma}, T, p, \|u_0\|_{\lambda_0,\infty}}  ,
\end{align}
%where the last constant $C$ is depend on $\lambda_0, c_1, L_b, K_{\sigma}, T, p, \|u_0\|_{\lambda_0,\infty}$,
where $\|u_0\|_{\lambda_0,\infty}:= \sup_{x\in\mathbb{R}}(|u_0(x)|e^{-\lambda_0 |x|})$. Note that the last constant in (\ref{210225.1854}) is independent of $n$.
Hence (\ref{210225.1905}) is proved, completing the proof of the lemma.
%Now, (\ref{05142144}) follows from (\ref{210124.2114}).  $\blacksquare$

\vskip 0.6cm

We will apply a Kolmogorov type tightness criterion (see Lemma 6.3 of \cite{S}) to establish the tightness of the law of $\{u_n\}$ in $C(\mathbb{R}_+, C_{tem})$.  This is given in the following lemma.

%Next, we use the uniform estimate (\ref{05142144}) to establish the tightness of the law of $\{u_n\}$ in $C(\mathbb{R}_+, C_{tem})$. To do this, we will apply a Kolmogorov tightness criterion, see {}. Therefore, it is sufficient to prove the following estimate.

%
%Next, we use the uniform estimate (\ref{05142144}) to establish the tightness of the law of $\{u_n\}$ in $C(\mathbb{R}_+, C_{tem})$.
%In (\ref{05150908}), the terms
%\begin{align}
%  P_t u_0(x) + \int_0^t\int_{\mathbb{R}} p_{t-r}(x,z)\sigma W(dr,dz)
%\end{align}
%are independent of $n$. Hence to show the tightness of $\{u_n\}$, it suffices to show the sequence of processes
%\begin{align}
%  Z_n(t,x)=\int_0^t\int_{\mathbb{R}} p_{t-r}(x,z) b_n(u_n(r,z)) drdz ,\quad n\geq 1
%\end{align}
%is tight in $C(\mathbb{R}_+, C_{tem})$. To do this, we will apply a Kolmogorov tightness criterion, see {}. Therefore, it is sufficient to prove the following estimate.

Define
\begin{align}
  X_n(t,x):=\int_0^t\int_{\mathbb{R}} p_{t-r}(x,z) b_n(u_n(r,z)) \,\mathrm{d}r\mathrm{d}z ,\quad n\geq 1 .
\end{align}

\begin{lemma}\label{05152151}
Let  $u_0 \in C_{tem}$. Assume that (H1) holds and that $\sigma$ is continuous with $K_{\sigma}:=\sup_{z\in\mathbb{R}}|\sigma(z)| < \infty$.
Then for any $\lambda, T>0$, $p\geq 1$ and $\theta\in(0,1)$, there exist constants $C_{\lambda, c_1, L_b, K_{\sigma},  T, p, \theta, u_0}$ and $C_{K_{\sigma}, T, p}$ independent of $n$ such that
\begin{align}
\label{05152147}
\mathbb{E}\left(|X_n(t,x)-X_n(s,y)|^p e^{-\lambda |x|}\right)\leq & C_{\lambda, c_1, L_b, K_{\sigma}, T, p, \theta, u_0}\left(|t-s|^{\theta p}+|x-y|^p\right), \\
\label{210125.2123}
\mathbb{E}\left(|V_n(t,x)-V_n(s,y)|^p e^{-\lambda |x|}\right)\leq & C_{K_{\sigma}, T, p}\left(|t-s|^{\frac{p}{4}}+|x-y|^{\frac{p}{2}}\right) ,
\end{align}
for any $s,t\in [0,T]$ and  $x,y\in\mathbb{R}$ with $|x-y|\leq 1$. In  particular,  the family $\{u_n\}$ is tight in $C(\mathbb{R}_+, C_{tem})$.
\end{lemma}
\noindent {\bf Proof}. It suffices to prove this lemma for sufficiently large $p$ and sufficiently small $\lambda$.
Fix $T>0$.  As $\lim_{\lambda\rightarrow 0}T^*(\lambda, c_1)=\infty$, there exists a positive constant $\lambda_T$  such that $T\leq T^*(\lambda, c_1)$ for all $\lambda\leq \lambda_T$. For any fixed $\lambda\leq\lambda_T$, we choose $\lambda_0 >0$ so that
\begin{align}\label{210225.1936}
	2\lambda_0 e^{\beta(\lambda_0, c_1) T} = \lambda .
\end{align}
In the following we write $\beta$ for $\beta(\lambda_0, c_1)$ to simplify the notation.
Without loss of generality, we assume $t\geq s$ .
{\allowdisplaybreaks\begin{align}\label{05152143}
  & |X_n(t,x)-X_n(s,y)|e^{-\lambda |x| } \nonumber\\
  \leq  & |X_n(t,x)-X_n(s,x)|e^{-\lambda |x| } + |X_n(s,x)-X_n(s,y)|e^{-\lambda |x| }\nonumber\\
  \leq & \left|\int_0^s\int_{\mathbb{R}} [p_{t-r}(x,z)-p_{s-r}(x,z)]b_n(u_n(r,z))\,\mathrm{d}r\mathrm{d}z\right|e^{-\lambda |x| } \nonumber\\
  & + \left|\int_s^t\int_{\mathbb{R}} p_{t-r}(x,z) b_n(u_n(r,z))\,\mathrm{d}r\mathrm{d}z\right|e^{-\lambda |x| } \nonumber\\
  & + \left|\int_0^s\int_{\mathbb{R}} [p_{s-r}(x,z)-p_{s-r}(y,z)] b_n(u_n(r,z))\,\mathrm{d}r\mathrm{d}z\right|e^{-\lambda |x| } \nonumber\\
  =: & \, J_1 +J_2 +J_3 .
\end{align}}
%Note that for any $\theta\in[0,1]$, $0\leq s\leq t$,
%\begin{align}
%  |p_t(z)-p_s(z)|\leq \frac{c^\theta |t-s|^{\theta}}{s^{\theta}}(p_s(z)+p_t(z)+p_{2t}(z))
%\end{align}

\noindent By (i) of Lemma \ref{A.4}, (\ref{05161718}) and (\ref{05132049.1}), we have
\begin{align}
  J_1
\leq & \int_0^s\int_{\mathbb{R}} |p_{t-r}(x,z)-p_{s-r}(x,z)|\times \nonumber\\
  & \left[ c_1 |u_n(r,z)|\log_+ |u_n(r,z)|+ L_b(|u_n(r,z)|+1) \right]\,\mathrm{d}r\mathrm{d}z \times e^{-\lambda |x|} \nonumber\\
\leq &  \int_0^s\int_{\mathbb{R}} \frac{(2\sqrt{2})^\theta |t-s|^{\theta}}{(s-r)^{\theta}}\big(p_{s-r}(x,z)+p_{t-r}(x,z)+p_{2(t-r)}(x,z)\big) \nonumber\\
  & \times \bigg\{ c_ 1 e^{\lambda_0 |z|e^{\beta r}}\lambda_0 |z|e^{\beta r}\times \sup_{z\in\mathbb{R}}\left(|u_n(r,z)|e^{-\lambda_0 |z|e^{\beta r}}\right) \nonumber\\
  & + c_1 e^{\lambda_0 |z|e^{\beta r}}\times\sup_{z\in\mathbb{R}}\left[|u_n(r,z)|e^{-\lambda_0 |z|e^{\beta r}}\log_+\left(|u_n(r,z)|e^{-\lambda_0 |z|e^{\beta r}}\right)\right] \nonumber\\
  & + L_b e^{\lambda_0 |z|e^{\beta r}}\times \sup_{z\in\mathbb{R}}\left(|u_n(r,z)|e^{-\lambda_0 |z|e^{\beta r}}\right) +L_b \bigg\} \,\mathrm{d}z\mathrm{d}r \times e^{-\lambda |x| } \nonumber\\
\leq & (2\sqrt{2})^{\theta}|t-s|^{\theta}\times \int_0^s \frac{dr}{(s-r)^{\theta}}\times  \bigg\{\sup_{r\leq T, z\in\mathbb{R}}\left(|u_n(r,z)|e^{-\lambda_0 |z|e^{\beta r}}\right) \nonumber\\
  & \times C_{\lambda_0,c_1,L_b,T}\left(e^{\lambda_0 |x|e^{\beta r}}\lambda_0 |x|e^{\beta r}+ e^{\lambda_0 |x|e^{\beta r}}\right)  \times e^{-\lambda|x| } \nonumber\\
  & + \sup_{r\leq T, z\in\mathbb{R}}\left[|u_n(r,z)|e^{-\lambda_0 |z|e^{\beta r}}\log_+\left(|u_n(r,z)|e^{-\lambda_0 |z|e^{\beta r}}\right)\right]\nonumber\\
  & \times C_{\lambda_0,c_1,L_b,T}\times e^{\lambda_0 |x|e^{\beta r}}\times e^{-\lambda|x|} + 3L_b \bigg\} \,\mathrm{d}r.
\end{align}
By the choice of $\lambda_0$,
\begin{align}\label{200405.2053}
  e^{\lambda_0 |x|e^{\beta r}}\lambda_0 |x|e^{\beta r}\times e^{-\lambda|x|} \leq 1, \quad \forall\, r\in[0,T].
\end{align}
Therefore
\begin{align}\label{4.1}
  J_1\leq & (2\sqrt{2})^{\theta}|t-s|^{\theta}\times \int_0^s \frac{dr}{(s-r)^{\theta}}\,\mathrm{d}r \times C_{\lambda_0,c_1,L_b,T}\nonumber\\
  &\times
  \Bigg\{\sup_{r\leq T, z\in\mathbb{R}}\left(|u_n(r,z)|e^{-\lambda_0 |z|e^{\beta r}}\right)  + \sup_{r\leq T, z\in\mathbb{R}}\left(|u_n(r,z)|e^{-\lambda_0 |z|e^{\beta r}}\right) \nonumber\\
  & \times \log_+\left[\sup_{r\leq T, z\in\mathbb{R}}\left(|u_n(r,z)|e^{-\lambda_0 |z|e^{\beta r}}\right)\right] + 1\Bigg\} .
\end{align}
By Lemma \ref{lemma 05151406}, (\ref{210225.1854}) and the fact that $\theta\in(0,1)$ , we deduce from (\ref{4.1}) that
\begin{align}\label{05152144.1}
  E J_1^p\leq C_{\lambda_0,c_1,L_b,K_{\sigma},T,p,\theta,\|u_0\|_{\lambda_0,\infty}} |t-s|^{\theta p}.
\end{align}
Similarly, we can show that
\begin{align}\label{05152144.2}
    E J_2^p\leq C_{\lambda_0,c_1,L_b,K_{\sigma},T,p,\|u_0\|_{\lambda_0,\infty}}|t-s|^{ p}.
\end{align}
Now we estimate $J_3$. By (ii)-(iv) of Lemma \ref{A.4}, we have
\begin{align}
  J_3
\leq &  \int_0^s\int_{\mathbb{R}} |p_{s-r}(x,z)-p_{s-r}(y,z)|\times \nonumber\\
  & \left[ c_1 |u_n(r,z)|\log_+ |u_n(r,z)|+ L_b(|u_n(r,z)|+1) \right] \,\mathrm{d}r\mathrm{d}z \times e^{-\lambda |x|} \nonumber\\
\leq &  \int_0^s\int_{\mathbb{R}} |p_{s-r}(x,z)-p_{s-r}(y,z)|  \nonumber\\
  & \times \bigg\{c_1 e^{\lambda_0 |z|e^{\beta r}}\lambda_0 |z|e^{\beta r}\times \sup_{z\in\mathbb{R}}\left(|u_n(r,z)|e^{-\lambda_0 |z|e^{\beta r}}\right) \nonumber\\
  & + c_1 e^{\lambda_0 |z|e^{\beta r}}\times\sup_{z\in\mathbb{R}}\left[|u_n(r,z)|e^{-\lambda_0 |z|e^{\beta r}}\log_{+}\left(|u_n(r,z)|e^{-\lambda_0 |z|e^{\beta r}}\right)\right] \nonumber\\
  & + L_b e^{\lambda_0 |z|e^{\beta r}}\times \sup_{z\in\mathbb{R}}\left(|u_n(r,z)|e^{-\lambda_0 |z|e^{\beta r}}\right) +L_b \bigg\} \,\mathrm{d}z\mathrm{d}r \times e^{-\lambda |x|} \nonumber\\
\leq & |x-y|\times\int_0^s \frac{1}{\sqrt{s-r}} \times  \bigg\{\sup_{r\leq T, z\in\mathbb{R}}\left(|u_n(r,z)|e^{-\lambda_0 |z|e^{\beta r}}\right)  \nonumber\\
  & \times C_{\lambda_0, c_1, L_b,T} \Big[e^{\lambda_0 (|x|+|x-y|)e^{\beta r}}\lambda_0 (|x|+|x-y|)e^{\beta r} + e^{\lambda_0 (|x|+|x-y|)e^{\beta r}}\Big]  \nonumber\\
  & + \sup_{r\leq T,z\in\mathbb{R}}\left(|u_n(r,z)|e^{-\lambda_0 |z|e^{\beta r}}\right)\times \log_{+}\left[\sup_{r\leq T,z\in\mathbb{R}}\left(|u_n(r,z)|e^{-\lambda_0 |z|e^{\beta r}}\right)\right]\nonumber\\
  & \times C_{\lambda_0,c_1,L_b,T}\times e^{\lambda_0(|x|+|x-y|)e^{\beta r}}  + \sqrt{2}L_b  \bigg\} \,\mathrm{d}r \times e^{-\lambda |x|} .
\end{align}
Due to the fact that $|x-y|\leq 1$ and the choice of $\lambda_0$, we have
\begin{align}\label{200405.20531}
  & e^{\lambda_0 (|x|+|x-y|) e^{\beta r}}\lambda_0 (|x|+|x-y|) e^{\beta r}\times e^{-\lambda|x|} \nonumber\\
  \leq & e^{\lambda_0 (|x|+1) e^{\beta r}}\lambda_0 (|x|+1) e^{\beta r}\times e^{-\lambda(|x|+1)} e^{\lambda} \nonumber\\
  \leq & e^{\lambda}, \quad \forall\, r\in[0,T].
\end{align}
Hence in view of (\ref{210225.1854}), we see that
\begin{align}\label{05152144.3}
  E J_3^p \leq C_{\lambda_0,c_1,L_b,K_{\sigma},T,p,\|u_0\|_{\lambda_0,\infty}} |x-y|^p .
\end{align}
Combining (\ref{05152143}), (\ref{05152144.1}), (\ref{05152144.2}) and (\ref{05152144.3}) together yields
\begin{align}
	\mathbb{E} \left( |X_n(t,x)-X_n(s,y)|e^{-\lambda |x|} \right)^p \leq C_{\lambda_0,c_1,L_b,K_{\sigma},T,p,\theta,\|u_0\|_{\lambda_0,\infty}}\left(|t-s|^{\theta p}+|x-y|^p\right)  ,
\end{align}
where the constant $\lambda_0$ is determined by $\lambda, T, c_1$ according to (\ref{210225.1936}). Thus,
 (\ref{05152147}) is proved.
%Taking $\lambda$ to be sufficiently small in the above inequality, we can see that  holds for any fixed $\lambda>0$.

Now we prove (\ref{210125.2123}). Observe that
\begin{align}
	& |V_n(t,x)-V_n(s,y)| \nonumber\\
	\leq & \bigg|\int_0^s\int_{\mathbb{R}}[p_{t-r}(x,z)-p_{s-r}(x,z)]\sigma_n(u_n(r,z))\,W(\mathrm{d}r,\mathrm{d}z)\bigg| \nonumber\\
	  &  + \bigg|\int_s^t\int_{\mathbb{R}}p_{t-r}(x,z)\sigma_n(u_n(r,z))\,W(\mathrm{d}r,\mathrm{d}z)\bigg| \nonumber\\
	  &  + \bigg|\int_0^t\int_{\mathbb{R}}[p_{s-r}(x,z)-p_{s-r}(y,z)]\sigma_n(u_n(r,z))\,W(\mathrm{d}r,\mathrm{d}z)\bigg| \nonumber\\
	=: & \, I_1 +I_2 + I_3 .
\end{align}
So
\begin{align}\label{210126.2129}
	|V_n(t,x)-V_n(s,y)|^p \leq 3^{p-1} (I_1^p + I_2^p + I_3^p) .
\end{align}
By the BDG inequality, (\ref{210122.1855}), and (v) of Lemma \ref{A.4}, we get
\begin{align}\label{210126.2130}
	\mathbb{E} I_1^p \leq & \mathbb{E}\bigg[\int_0^s\int_{\mathbb{R}}\big|p_{t-r}(x,z)-p_{s-r}(x,z)\big|^2\sigma_n(u_n(r,z))^2 \,\mathrm{d}z\mathrm{d}r\bigg]^{\frac{p}{2}} \nonumber\\
	\leq & K_{\sigma}^p \times\bigg[\int_0^s\int_{\mathbb{R}}\big|p_{t-r}(x,z)-p_{s-r}(x,z)\big|^2  \,\mathrm{d}z\mathrm{d}r\bigg]^{\frac{p}{2}} \nonumber\\
	\leq & C_{K_{\sigma},T,p} |t-s|^{\frac{p}{4}}.
\end{align}
%where we have used (\ref{210222.1639}).
%there exists a constant $C>0$ such that
%\begin{align}\label{210126.2120}
%	\int_0^{t\vee s}\int_{\mathbb{R}} |p_{t-r}(x,z)-p_{s-r}(y,z)|^2 dzdr \leq C(|t-s|^{\frac{1}{2}}+|x-y|) .
%\end{align}
Similarly, we have
\begin{align}\label{210126.2131}
	\mathbb{E}I_3^p \leq C_{K_{\sigma},T,p} |x-y|^{\frac{p}{2}} .
\end{align}
For the term $I_2$, the uniform boundedness of $\sigma_n$ gives
\begin{align}\label{210126.2132}
	\mathbb{E} I_2^p \leq & \mathbb{E}\bigg[\int_s^t\int_{\mathbb{R}} p_{t-r}(x,z)^2\sigma_n(u_n(r,z))^2 \,\mathrm{d}z\mathrm{d}r\bigg]^{\frac{p}{2}} \nonumber\\
	\leq & K_{\sigma}^p \times\bigg[\int_s^t\int_{\mathbb{R}}\big|p_{t-r}(x,z)\big|^2  \,\mathrm{d}z\mathrm{d}r\bigg]^{\frac{p}{2}} \nonumber\\
	\leq & K_{\sigma}^p \times\bigg[\int_s^t \frac{1}{2\sqrt{\pi(t-r)}}  \,\mathrm{d}r\bigg]^{\frac{p}{2}} \nonumber\\
	\leq & C_{K_{\sigma},T,p} |t-s|^{\frac{p}{4}} .
\end{align}
Combining (\ref{210126.2129}), (\ref{210126.2130}), (\ref{210126.2131}) and (\ref{210126.2132}) together, we obtain (\ref{210125.2123}) .
This completes the proof of Lemma \ref{05152151}.
$\blacksquare$

\vskip 0.6cm

%\begin{theorem}
%If $u_0\in C_{tem}$, then there exists a solution to (\ref{1.a}), and the sample paths of the solution are a.e. in $C(\mathbb{R}_{+}, C_{tem})$.
%\end{theorem}
\noindent {\bf Proof of Theorem \ref{thm1}}.
%Note that
%\begin{align}
%  u_n(t,x) = P_t u_0(x)  + X_n(t,x) + V_n(t,x) ,
%\end{align}
%where the first terms on the right hand side of the above equality is independent of $n$, $\{X_n\}$ and $\{V_n\}$ are tight in $C(\mathbb{R}_+,C_{tem})$ by Lemma \ref{05152151}. Therefore, we obtain the tightness of $\{u_n\}$ in $C(\mathbb{R}_+,C_{tem})$.
We have established in Lemma \ref{05152151} that $\{u_n\}$ is tight in $C(\mathbb{R}_+,C_{tem})$.
By Prokhorov's theorem and the modified version of Skorokhod's representation theorem whose proof can be found in Appendix C of [BHR], we may assume that $d(u_n, u)\rightarrow 0$ (not relabelled) a.s. in $C(\mathbb{R}_+,C_{tem})$ for some process $u$ on some probability space $(\widetilde{\Omega}, \widetilde{\mathcal{F}}, \widetilde{\mathbb{P}})$, in other words, for any $\lambda>0$, $T\geq 0$,
\begin{align}\label{05160930.2}
  \sup_{t\leq T, x\in\mathbb{R}}\left(|u_n(t,x)-u(t,x)|e^{-\lambda|x|}\right)\rightarrow 0, \quad \widetilde{\mathbb{P}}-a.s..
\end{align}
By the dominated convergence theorem, and using (\ref{05160930.2}), (\ref{05160930.1}), (\ref{210126.2140}),  (\ref{05160929.2}) and (\ref{210122.1855}),  one can deduce that for any $(t,x)\in\mathbb{R_+}\times\mathbb{R}$,
\begin{align}
  \int_0^t\int_{\mathbb{R}} p_{t-r}(x,z)b_n(u_n(r,z))\,\mathrm{d}r\mathrm{d}z\rightarrow \int_0^t\int_{\mathbb{R}} p_{t-r}(x,z)b(u(r,z))\,\mathrm{d}r\mathrm{d}z ,
\end{align}
$\widetilde{\mathbb{P}}$-a.s. as $n\rightarrow\infty$, and
\begin{align}
	\int_0^t\int_\mathbb{R} p_{t-r}(x,z)\sigma_n(u_n(r,z))\,W(\mathrm{d}r,\mathrm{d}z)\rightarrow \int_0^t\int_{\mathbb{R}} p_{t-r}(x,z)\sigma(u(r,z))\,W(\mathrm{d}r,\mathrm{d}z),
\end{align}
in the sense of $L^p(\widetilde{\Omega})$ for any $p\geq 1$ as $n\rightarrow\infty$.
Therefore, we see that $u$ is a  weak solution of (\ref{1.a}).
%The existence of strong solutions in the probabilistic sense follows from the path uniqueness of the solutions and the Yamada-Watanabe theorem.
$\blacksquare$

\section{Pathwise uniqueness}
In this section, we prove the pathwise uniqueness of solutions to (\ref{1.a}) and hence obtain the strong solution.
%We first show the pathwise uniqueness of the solutions in $C(\mathbb{R}_{+},C_{tem})$. Then

%\begin{lemma}\label{lemma 4.1}
%The pathwise uniqueness holds for the solutions of  (\ref{1.a}) in $C(\mathbb{R}_{+}, C_{tem})$.
%\end{lemma}

\vskip 0.3cm
\noindent {\bf Proof of Theorem \ref{thm3}}. Since condition (H2) implies condition (H1), there exists a weak solution to (\ref{1.a}) according to Theorem \ref{thm1}. We only  show the pathwise uniqueness for solutions of (\ref{1.a}). The existence of strong solutions then follows from the Yamada-Watanabe theorem.

Suppose that $u,v$ are two solutions of equation (\ref{1.a}) that belong to the space $C(\mathbb{R}_{+}, C_{tem})$. We are going to show that $u=v$.
%To prove pathwise uniqueness, it suffices to show that for any fixed $T>0$, $u(t,x)=v(t,x)$ for any $x\in\mathbb{R}$ and any $$
Fix $T>0$, and take $\lambda>0$ sufficiently small so that $T\leq T^*(\lambda, c_4)$. In this section, we write $\beta$ for $\beta(\lambda, c_4)$ for simplicity.
%To prove the pathwise uniqueness, Suppose that $u,v$ are any two solutions of (\ref{1.a}) and their paths are in $C(\mathbb{R}_{+}, C_{tem})$.
%Let $\delta$ be the constant appeared in Condition (H2).
%Without loss of generality, we can assume $0<\delta\leq e^{-1}$.
Let $M>0 $ and $0<\delta\leq e^{-1}$.
Define stopping times
\begin{align*}
  \tau_M := & \inf\left\{t>0: \sup_{x\in\mathbb{R}}\left(|u(t,x)|e^{-\lambda|x|e^{\beta t}}\right)\geq M\right\} \\
  &\wedge \inf\left\{t>0: \sup_{x\in\mathbb{R}}\left(|v(t,x)|e^{-\lambda|x|e^{\beta t}}\right)\geq M\right\} , \\
  \tau^{\delta}:=& \inf\left\{t>0: \sup_{x\in\mathbb{R}}\left(|u(t,x)-v(t,x)|e^{-\lambda|x|e^{\beta t}}\right)\geq \delta\right\} , \\
  \tau_M^{\delta}:=& \tau_M\wedge\tau^{\delta}\wedge T ,
\end{align*}
with the convention that $\inf \emptyset = +\infty$.
Define also
\begin{align}
  Z(r):= \mathbb{E}\sup_{t\leq r\wedge\tau_M^{\delta}, x\in\mathbb{R}}\left(|u(t,x)-v(t,x)|e^{-\lambda |x|e^{\beta t}}\right) .
\end{align}

\noindent Obviously
%{\allowdisplaybreaks \begin{align}
%  & u(t,x)-v(t,x) \nonumber\\
%  =& \int_0^t\int_{\mathbb{R}} p_{t-s}(x,y) e^{\lambda |y| e^{\beta s}}\lambda |y| e^{\beta s}\left(u(s,y)e^{-\lambda |y| e^{\beta s}}\right) dsdy \nonumber\\
%  & + \int_0^t\int_{\mathbb{R}} p_{t-s}(x,y) e^{\lambda |y| e^{\beta s}}\left(u(s,y)e^{-\lambda |y| e^{\beta s}}\right)\log\left(|u(s,y)|e^{-\lambda |y| e^{\beta s}}\right) dsdy \nonumber\\
%  & + \int_0^t\int_{\mathbb{R}} p_{t-s}(x,y) \sigma(u(s,y)) W(ds,dy) \nonumber\\
%  & - \int_0^t\int_{\mathbb{R}} p_{t-s}(x,y) e^{\lambda |y| e^{\beta s}}\lambda |y| e^{\beta s}\left(v(s,y)e^{-\lambda |y| e^{\beta s}}\right) dsdy \nonumber\\
%  & - \int_0^t\int_{\mathbb{R}} p_{t-s}(x,y) e^{\lambda |y| e^{\beta s}}\left(v(s,y)e^{-\lambda |y| e^{\beta s}}\right)\log\left(|v(s,y)|e^{-\lambda |y| e^{\beta s}}\right) dsdy \nonumber\\
%  & - \int_0^t\int_{\mathbb{R}} p_{t-s}(x,y) \sigma(v(s,y)) W(ds,dy) .
%\end{align}
%}
%{\allowdisplaybreaks \begin{align}
%  & u(t,x)-v(t,x) \nonumber\\
%  =& \int_0^t\int_{\mathbb{R}} p_{t-s}(x,y) [b(u(s,y))-b(v(s,y))] dsdy \nonumber\\
%  & + \int_0^t\int_{\mathbb{R}} p_{t-s}(x,y) [\sigma(u(s,y))-\sigma(v(s,y))] W(ds,dy) \nonumber\\
%  =: & I + J .
%\end{align}
%}
%{\allowdisplaybreaks
\begin{align}\label{210224.1634}
  & Z(r) \nonumber\\
  \leq & \,\mathbb{E}\sup_{t\leq r\wedge\tau_M^{\delta}, x\in\mathbb{R}} \left\{  \int_0^t\int_{\mathbb{R}} p_{t-s}(x,y)  |b(u(s,y))-b(v(s,y)) | \,\mathrm{d}s\mathrm{d}y \times e^{-\lambda |x|e^{\beta t}} \right\}\nonumber\\
  & + \mathbb{E}\sup_{t\leq r\wedge\tau_M^{\delta}, x\in\mathbb{R}} \left\{ \left| \int_0^t\int_{\mathbb{R}} p_{t-s}(x,y) [\sigma(u(s,y))-\sigma(v(s,y))] \,W(\mathrm{d}s,\mathrm{d}y) \right| e^{-\lambda |x|e^{\beta t}} \right\}  \nonumber\\
  =: & \,  I + J .
\end{align}
%}
Now we estimate the term $I, J$ separately. By condition (H2), we have
{\allowdisplaybreaks\begin{align}\label{210224.1635}
	I \leq & \,\mathbb{E}\sup_{t\leq r\wedge\tau_M^{\delta}, x\in\mathbb{R}}\bigg\{\int_0^t \int_{\mathbb{R}} p_{t-s}(x,y) c_3|u(s,y)-v(s,y)| \nonumber\\
	& ~~~~~~~\times \log_{+}\frac{1}{|u(s,y)-v(s,y)|} \,\mathrm{d}s\mathrm{d}y \times  e^{-\lambda |x| e^{\beta t}} \bigg\} \nonumber\\
	& +  \mathbb{E}\sup_{t\leq r\wedge\tau_M^{\delta}, x\in\mathbb{R}}\bigg\{\int_0^t \int_{\mathbb{R}} p_{t-s}(x,y) c_4 \log_{+}\big(|u(s,y)|\vee |v(s,y)|\big) \nonumber\\
	& ~~~~~~~ \times |u(s,y)-v(s,y)| \,\mathrm{d}s\mathrm{d}y \times  e^{-\lambda |x| e^{\beta t}} \bigg\} \nonumber\\
	& + \mathbb{E}\sup_{t\leq r\wedge\tau_M^{\delta}, x\in\mathbb{R}}\bigg\{\int_0^t \int_{\mathbb{R}} p_{t-s}(x,y) c_5 |u(s,y)-v(s,y)| \,\mathrm{d}s\mathrm{d}y \times  e^{-\lambda |x| e^{\beta t}} \bigg\} \nonumber\\
	=: & \, I_1 + I_2 + I_3 .
\end{align}}

\noindent First, we estimate the term $I_1$. By the fact that the function $x\mapsto x\log\frac{1}{x}$ is increasing and concave on $(0,e^{-1})$, (\ref{05132049.1}) and (\ref{05132049.2}), we get
{\allowdisplaybreaks\begin{align}
	I_1 \leq & c_3 \,\mathbb{E}\sup_{t\leq r\wedge\tau_M^{\delta}, x\in\mathbb{R}}\bigg\{\int_0^t \int_{\mathbb{R}} p_{t-s}(x,y) e^{\lambda |y| e^{\beta s}} |u(s,y)-v(s,y)| e^{-\lambda |y| e^{\beta s}} \nonumber\\
	& ~~~~~~~\times \log_{+}\frac{1}{|u(s,y)-v(s,y)| e^{-\lambda |y| e^{\beta s}}} \,\mathrm{d}y\mathrm{d}s \times  e^{-\lambda |x| e^{\beta t}} \bigg\} \nonumber\\
	\leq & c_3  \,\mathbb{E}\sup_{t\leq r\wedge\tau_M^{\delta}, x\in\mathbb{R}}\bigg\{\int_0^t \sup_{y\in\mathbb{R}} \bigg[|u(s,y)-v(s,y)| e^{-\lambda |y| e^{\beta s}} \nonumber\\
	& ~~~~~~~ \times \log_{+}\frac{1}{|u(s,y)-v(s,y)| e^{-\lambda |y| e^{\beta s}}}\bigg] \nonumber\\
	& ~~~~~~~ \times \int_{\mathbb{R}} p_{t-s}(x,y) e^{\lambda |y| e^{\beta s}} \,\mathrm{d}y\mathrm{d}s \times  e^{-\lambda |x| e^{\beta t}} \bigg\} \nonumber\\
	\leq & 2c_3 e^{\frac{\lambda^2}{4\beta} e^{2\beta r-1}} \int_0^t \mathbb{E} \bigg\{ \sup_{\rho\leq s\wedge\tau_M^{\delta}, y\in\mathbb{R}} \Big( |u(\rho,y)-v(\rho,y)| e^{-\lambda |y| e^{\beta \rho}} \Big) \nonumber\\
	& \times \log_{+} \frac{1}{\sup_{\rho\leq s\wedge\tau_M^{\delta}, y\in\mathbb{R}}\big(|u(\rho,y)-v(\rho,y)| e^{-\lambda |y| e^{\beta \rho}}\big)} \bigg\} \,\mathrm{d}s \nonumber\\
	\leq & 2c_3 e^{\frac{\lambda^2}{4\beta} e^{2\beta r-1}} \int_0^r Z(s)\log_{+}\frac{1}{Z(s)} \,\mathrm{d}s,
\end{align}}

\noindent where (\ref{05132049.1}) was used.
Note that
\begin{align*}
	\log_{+}(ab) \leq \log_{+}a + \log_{+} b .
\end{align*}
By the definition of $\tau_M^{\delta}$, we have
{\allowdisplaybreaks\begin{align}
	I_2
%	\leq & \mathbb{E}\sup_{t\leq r\wedge\tau_M^{\delta}, x\in\mathbb{R}}\bigg\{\int_0^t \int_{\mathbb{R}} p_{t-s}(x,y) c_4 \log_{+}\big(|u(s,y)|\vee |v(s,y)|\big) \nonumber\\
%	& ~~~~~~~ \times |u(s,y)-v(s,y)| dsdy \times  e^{-\lambda |x| e^{\beta t}} \bigg\} \nonumber\\
	\leq & c_4 \,\mathbb{E}\sup_{t\leq r\wedge\tau_M^{\delta}, x\in\mathbb{R}}\bigg\{\int_0^t \int_{\mathbb{R}} p_{t-s}(x,y) \bigg[ \log_{+}\Big(e^{\lambda|y| e^{\beta s} }\Big) \nonumber\\
	& + \log_{+}\Big(\big(|u(s,y)|e^{-\lambda|y| e^{\beta s}}\big)\vee \big(|v(s,y)| e^{-\lambda|y| e^{\beta s}}\big) \Big)  \bigg]\nonumber\\
	&  \times \Big(|u(s,y)-v(s,y)|e^{-\lambda|y| e^{\beta s}}\Big) e^{\lambda|y| e^{\beta s}} \,\mathrm{d}y\mathrm{d}s \times  e^{-\lambda |x| e^{\beta t}} \bigg\} \nonumber\\
	\leq & c_4 \,\mathbb{E}\sup_{t\leq r\wedge\tau_M^{\delta}, x\in\mathbb{R}}\bigg\{\int_0^t \sup_{y\in\mathbb{R}} \Big(|u(s,y)-v(s,y)|e^{-\lambda|y| e^{\beta s}}\Big)\nonumber\\
	& \times  \int_{\mathbb{R}} p_{t-s}(x,y) e^{\lambda|y| e^{\beta s}} \lambda|y| e^{\beta s} \,\mathrm{d}y\mathrm{d}s  \times  e^{-\lambda |x| e^{\beta t}} \bigg\} \nonumber\\
	& + c_4 \log_{+}(M)  \mathbb{E}\sup_{t\leq r\wedge\tau_M^{\delta}, x\in\mathbb{R}}\bigg\{\int_0^t \sup_{y\in\mathbb{R}} \Big(|u(s,y)-v(s,y)|e^{-\lambda|y| e^{\beta s}}\Big)\nonumber\\
	& \times  \int_{\mathbb{R}} p_{t-s}(x,y) e^{\lambda|y| e^{\beta s}} \,\mathrm{d}y\mathrm{d}s  \times  e^{-\lambda |x| e^{\beta t}} \bigg\} \nonumber\\
	\leq & \frac{1}{2} Z(r) + c_4 C_{\lambda, \beta, M, r} \int_0^r Z(s)\,\mathrm{d}s ,
\end{align}}

\noindent where the last inequality holds for the same reason as the derivation of (\ref{05132036})-(\ref{05132040.2}) with constant $c_1$  replaced by constant $c_4$.
Similarly,
\begin{align}
	I_3 \leq 2 c_5 e^{\frac{\lambda^2}{4\beta} e^{2\beta r-1}} \int_0^r Z(s) \,\mathrm{d}s .
\end{align}

\noindent For the term $J$, we use the estimate established in Proposition \ref{estimates 003} to obtain
\begin{align*}
	J \leq & \epsilon \,\mathbb{E} \sup_{s\leq r\wedge\tau_M^{\delta}, y\in\mathbb{R}}\left(|\sigma(u(s,y))-\sigma(v(s,y))| e^{-\lambda|y|e^{\beta s}}\right) \nonumber\\
	& + C_{\epsilon,\lambda, \beta, r} \,\mathbb{E} \int_0^{r\wedge\tau_{M}^{\delta}}\int_{\mathbb{R}}|\sigma(u(s,y))-\sigma(v(s,y))|  e^{-\lambda |y|e^{\beta s}} \,\mathrm{d}y\mathrm{d}s  ,
\end{align*}
where the constant $C_{\epsilon,\lambda, \beta, r}$ is the constant $C_{\epsilon,p,h(T),T}$ appeared in (\ref{101.2}) by taking $p=1$, $T=r$ and $h(T)=\lambda e^{\beta r}$.
Since $\sigma$ is bounded and Lipschitz, there exists two nonnegative constants $K_{\sigma}$ and $L_{\sigma}$ such that
\begin{align*}
	|\sigma(x)| \leq & K_{\sigma}, \quad \forall\, x\in\mathbb{R}, \\
	|\sigma(x)-\sigma(y)| \leq & L_{\sigma}|x-y|, \quad \forall\, x,y\in\mathbb{R}.
\end{align*}
Hence for any $0< \theta <1$, we have
\begin{align}\label{210217.2114}
	J \leq & \epsilon L_{\sigma} Z(r) + C_{\epsilon,\lambda, \beta, r}\,\mathbb{E} \int_0^{r\wedge\tau_{M}^{\delta}} \sup_{y\in\mathbb{R}} \bigg\{ \left(|\sigma(u(s,y))-\sigma(v(s,y))| e^{-\lambda |y|e^{\beta s}} \right)^{\theta} \nonumber\\
	& \times \int_{\mathbb{R}} |\sigma(u(s,y))-\sigma(v(s,y))|^{1-\theta}  e^{-(1-\theta)\lambda |y|e^{\beta s}} \,\mathrm{d}y \bigg\} \,\mathrm{d}s \nonumber\\
	\leq & \epsilon L_{\sigma} Z(r) + \frac{(2K_{\sigma})^{1-\theta} L_{\sigma}^{\theta} C_{\epsilon,\lambda, \beta, r}}{(1-\theta)\lambda} \mathbb{E} \int_0^{r\wedge\tau_{M}^{\delta}} \sup_{y\in\mathbb{R}} \left(| u(s,y)- v(s,y)| e^{-\lambda |y|e^{\beta s}} \right)^{\theta} \mathrm{d}s \nonumber\\
	\leq & \epsilon L_{\sigma} Z(r) + \frac{(2K_{\sigma})^{1-\theta} L_{\sigma}^{\theta} C_{\epsilon,\lambda, \beta, r}}{(1-\theta)\lambda} \int_0^r Z(s)^{\theta} \,\mathrm{d}s .
\end{align}
Combining (\ref{210224.1634})-(\ref{210217.2114}) together, we obtain that
\begin{align}
  Z(r)\leq & \left(\frac{1}{2} + \epsilon L_{\sigma} \right) Z(r) +  C_{\lambda,M, c_4, c_5, r} \int_0^r Z(s)\,\mathrm{d}s \nonumber\\
  & + 2c_3 e^{\frac{\lambda^2}{4\beta}e^{2\beta r-1}} \int_0^r Z(s)\log_{+}\frac{1}{Z(s)}\,\mathrm{d}s +\frac{(2K_{\sigma})^{1-\theta} L_{\sigma}^{\theta} C_{\epsilon,\lambda, \beta, r}}{(1-\theta)\lambda} \int_0^r Z(s)^{\theta} \,\mathrm{d}s .
\end{align}
Taking for example $\epsilon= \frac{1}{4L_{\sigma}}$, subtracting $\left(\frac{1}{2} + \epsilon L_{\sigma} \right) Z(r)$ from both sides of the above inequality, and then applying the special Gronwall-type inequality established in Lemma \ref{A.3}, we obtain
\begin{align}
  Z(r)\equiv 0, \quad \forall\, r\geq 0 .
\end{align}
Since the solutions of (\ref{1.a}) don't blowup, let $M\rightarrow\infty$ to obtain  $\mathbb{P}$-a.s.,
\begin{align}
  u(t,x)=v(t,x), \quad \forall\, x\in \mathbb{R}, \ \forall\, t\in [0, \tau^{\delta}\wedge T].
\end{align}
%$u(t,x)=v(t,x)$ for any $x\in \mathbb{R}$ and any $t\in [0, \tau^{\delta}]$.
This implies that $\tau^{\delta} \geq T $, $\mathbb{P}$-a.s., otherwise it contradicts the definition of $\tau^{\delta}$. By the arbitrariness of $T$, we obtain that for $\mathbb{P}$-a.s.,
\begin{align}
	u(t,x)=v(t,x), \quad \forall\, (t,x)\in \mathbb{R}_{+}\times\mathbb{R}.
\end{align}
This completes the proof the pathwise uniqueness.
$\blacksquare$

\vskip 0.6cm

\noindent{\bf Acknowledgement}. We are grateful to Robert Dalang and Davar Khoshnevisan for their useful suggestions and comments. This work is partially supported by by NSFC (No. 11971456, 11721101, No. 12001516).

\end{document}